\documentclass[11pt]{amsart}  
\usepackage[T1]{fontenc}
\usepackage{amsxtra}

\usepackage{color}
\usepackage{amscd,amsfonts,amsmath,amssymb,
amsthm,amstext,latexsym,amstext}
\usepackage{enumerate,pifont,graphicx}
\usepackage[english]{babel}
\usepackage[all,cmtip]{xy} % for commutative diagram
\newcommand{\Z}{{\mathbb Z}}

\newcommand{\R}{{\mathbb R}}
\newcommand{\C}{{\mathbb C}}

\renewcommand{\S}{{\mathbb S}}

\def\cqfd{\hfill$\Box$}
\def\Res{{\,\rm Res}}

\def\Re{{\rm Re}}
\def\Im{{\rm Im}}
\def\whu{{\widehat{u}}}
\def\whOmega{{\widehat{\Omega}}}

\def\wtu{{\widetilde{u}}}
\def\wtU{{\widetilde{U}}}
\def\wtOmega{{\widetilde{\Omega}}}

\newtheorem{theorem}{Theorem}
\newtheorem{lemma}{Lemma}
\newtheorem{proposition}{Proposition}
\newtheorem{remark}{Remark}

\newtheorem{claim}{Claim}

\newtheorem{definition}{Definition}

\title{Classification of the solutions to an overdetermined elliptic problem in the plane}
\renewcommand\shorttitle{Classification of the solutions to an overdetermined elliptic problem}
\author{Martin Traizet}
\begin{document}
\maketitle
\begin{center}
March 22, 2013
\end{center}

\bigskip

\bigskip

{\em Abstract: we classify the solutions to an overdetermined elliptic problem in
the plane in the finite connectivity case.
This is achieved by establishing a one-to-one correspondence between
the solutions to this problem and a certain type of minimal surfaces.
}
\medskip

\section{Introduction}
\label{section1}
In the theory of elliptic P.D.E., an overdetermined problem is one where both the
Dirichlet and Neumann boundary values are prescribed.
This puts strong geometric constraints on the domain.
For example, a famous result of J. Serrin \cite{serrin} asserts that if $\Omega$ is a
bounded domain in $\R^n$ which admits a function $u$ solution of $\Delta u=-1$
in $\Omega$ with Dirichlet boundary value $u=0$ on $\partial\Omega$
and Neumann boundary value
$\frac{\partial u}{\partial\nu}$ constant on $\partial\Omega$, then $\Omega$ is a ball.
Overdetermined elliptic problems appear in various mathematical and physical problems,
such as isoperimetric inequalities, spectral geometry and hydrodynamics. (See for example the survey \cite{beneteau}.)

\medskip

In \cite{hhp}, F. H\'elein, L. Hauswirth and F. Pacard have proposed the following overdetermined problem.
Let $\Omega$ be a smooth, unbounded domain in $\R^n$ with non-empty boundary.
The domain $\Omega$ is called exceptional if it admits a positive harmonic function
$u$ which has Dirichlet boundary value $u=0$ on $\partial\Omega$ and
Neumann boundary value $\frac{\partial u}{\partial\nu}=c$ on $\partial\Omega$, where
$c$ is a constant.
By the boundary maximum principle, the constant $c$ must be negative,
and we may normalize $c=-1$ by scaling of $u$.
Also, when $u$ is constant on $\partial\Omega$, the Neumann boudary condition is equivalent to $|\nabla u|=|c|$ on $\partial\Omega$.
So we may formulate the above problem as:
\begin{equation}
\label{eq1}
\left\{\begin{array}{ll}
\Delta u=0& \mbox{ in $\Omega$,}\\
u>0 &\mbox{ in $\Omega$,}\\
u=0 &\mbox{ on $\partial\Omega$,}\\
|\nabla u|=1& \mbox{ on $\partial\Omega$}.
\end{array}\right.\end{equation}
This problem is related to the study of extremal domains, namely domains $\Omega$
in a Riemannian manifold
which are critical points for the functional $\lambda_1(\Omega)$ under a volume constraint,
where $\lambda_1$ denotes the first eigenvalue of the Laplace Beltrami operator.
See \cite{hhp} for more details.

\medskip
For example, a half-space or the complementary of a ball are exceptional domains.
In \cite{hhp}, the authors discovered that in the case $n=2$,
the domain $|y|<\frac{\pi}{2}+\cosh x$
in the plane is an exceptional domain. They also developed a Weierstrass type
Representation for exceptional domains in the plane that are simply connected, and noted a strong
analogy with minimal surfaces.
In this paper, we prove that the analogy goes very deep by establishing a one-to-one
correspondence between exceptional domains
and a certain type of minimal surfaces which we call minimal bigraphs.
This correspondence allows us to find examples and classify solutions.

\medskip

We only address Problem \eqref{eq1} in the case of planar domains ($n=2$)
and we identify $\R^2$ with the complex plane $\C$.
The assumption that $\Omega$ is a smooth domain can be relaxed
(see Proposition \ref{proposition-smooth}).
Also, the solution $u$ of Problem \eqref{eq1}, when it exists, is unique (see Proposition \ref{proposition-unique}).

\medskip

By a trivial exceptional domain, we mean a half-plane.
Let $\Omega$ be a non-trivial exceptional domain.
In Section \ref{section2}, we use the theory of univalent functions to prove that

\medskip

{\em
\begin{itemize}
\item If $\Omega$ has finite connectivity, then $|\nabla u|< 1$ in $\Omega$.
\item If $\Omega$ is periodic and has finite connectivity in the quotient,
then $|\nabla u|<1$ in $\Omega$.
\end{itemize}}

\medskip

By {\em finite connectivity}, we mean that $\partial\Omega$ has a finite number of
components. By {\em periodic}, we mean that $\Omega$ is invariant by a non-zero translation $T$. Note that an exceptional domain cannot be doubly periodic
(for the maximum principle implies that $u\equiv 0$ in this case).
\medskip

In Sections \ref{section3}, \ref{section4} and \ref{section4bis}, we establish a one-to-one correspondence between the following two classes of objects:

\medskip

{\em 
\begin{itemize}
\item exceptional domains  $\Omega$ such that  $|\nabla u|<1$ in 
$\Omega$,
\item complete, embedded minimal surfaces $M$ in $\R^3$ which are symmetric
with respect to the horizontal plane $x_3=0$ and such that $M^+=M\cap \{x_3>0\}$
is a graph over the unbounded domain in the plane bounded by $M\cap\{x_3=0\}$.
We call such a minimal surface a minimal bigraph.
\end{itemize}
}

\medskip

(In fact, we will establish the above correspondence assuming that
the domain $\Omega$ satisfies a mild additional geometric hypothesis, namely
that its complement in non-thinning: see Definition \ref{definition-thin}.
This hypothesis is always satisfied for domains with finite connectivity, or
periodic domains with finite connectivity in the quotient).
For example:
\begin{enumerate}
\item The vertical catenoid is a minimal bigraph. It corresponds to the exceptional domain  $\Omega=
\C\setminus D(0,1)$.
\item The horizontal catenoid is a minimal bigraph. It corresponds to the
exceptional domain $|y|<\frac{\pi}{2}+\cosh x$ discovered in \cite{hhp}.
\item Scherk's family of simply periodic minimal surfaces, suitably rotated, are minimal bigraphs. They correspond to a new family of periodic exceptional domains.
(In fact, this family was already known! See Section 
\ref{section-related}.)
\end{enumerate}
We give more details about these examples in Section \ref{section5}.

\medskip

In Section \ref{section6}, we take advantage of the correspondence to translate 
classification results from
minimal surface theory into classification results for exceptional domains.
We prove that (up to similitude)
{\em
\begin{itemize}
\item The only exceptional domains in the plane with finite connectivity are
the half-plane and Examples (1) and (2),
\item The only periodic exceptional domains with finite connectivity in the quotient
are the half-plane and Examples (3).
\end{itemize}
}
\medskip

Finally, in Section \ref{section7}, we extend the correspondence to the case of
immersed domains in the plane.

\subsection{Related works}

\medskip

\label{section-related}

Laurent Hauswirth pointed out to me that Problem \eqref{eq1} has been studied by 
D. Khavinson, E. Lundberg and R. Teodorescu in a recent paper
\cite{khavinson}.
They obtain classification results in the 
$2$-dimensional case under
stronger topological hypotheses than ours. They prove that if an exceptional
domain $\Omega$ is the complement of a bounded domain, then it is
the complement of a disk; if it is simply-connected and Smirnov, then it is a half-plane
or the domain $|y|<\frac{\pi}{2}+\cosh x$, up to similitude.
In the simply-connected case, their results are
stronger than ours because they do not assume that the boundary of $\Omega$ has
a finite number of components.
They also prove that in higher dimension,
an exceptional domain in $\R^n$ whose complement is bounded, connected and
has $C^{2,\alpha}$ boundary, is the exterior of a sphere.
\medskip

I learned from Erik Lundberg that the family of periodic exceptional domains
(the ones corresponding to Scherk surfaces) already appears in a 1976 paper by
G. R. Baker, P. G. Saffman and J. S. Sheffield \cite{baker} as a solution to an equilibrium problem in hydrodynamics. It is also discussed by D. Crowdy and C. Green in \cite{crowdy}.

\medskip

Several authors have studied the following, more general overdetermined elliptic problem
\begin{equation}
\label{eq2}
\left\{\begin{array}{ll}
\Delta u =f(u) & \mbox{ in $\Omega$}\\
u>0 & \mbox{ in $\Omega$}\\
u=0 & \mbox{ on $\partial\Omega$}\\
\frac{\partial u}{\partial \nu}=c & \mbox{ on $\partial\Omega$}
\end{array}\right.
\end{equation}
where $\Omega$ is a domain in $\R^n$ and $f$ is some given function.
A formal analogy between this kind of problem and constant mean curvature (CMC) hypersurfaces
in $\R^{n+1}$ has been observed.
For example, the quoted result of J. Serrin \cite{serrin} is the counterpart of the
theorem of A.D. Alexandrov, which asserts that the only embedded compact CMC hypersurfaces
in $\R^{n+1}$ are round spheres.
In \cite{schlenk-sicbaldi}, F. Schlenk and P. Sicbaldi have constructed solutions to Problem \eqref{eq2}
in the case $f(x)=\lambda x$, which are analogues of the Delaunay CMC surfaces
in $\R^3$.
In the case $n=2$, the analogy between Problem \eqref{eq2} and CMC surfaces
has been explored in a systematic way by A. Ros and P. Sicbaldi in a very interesting
paper \cite{ros-sicbaldi}. See the first section of this paper for more results in this
spirit and related conjectures.
\section{Preliminary remarks}
\label{section-preliminaires}
The assumption that $\Omega$ is a smooth domain can be relaxed.
Recall (\cite{GT} page 94)
that an open set $\Omega$ in $\R^n$ with non-empty boundary
is a domain of class $C^k$ (resp. smooth, analytic...) if for each point $x_0\in\partial\Omega$,
there exists $\varepsilon>0$ and a diffeomorphisme $\psi$ of class $C^k$ (resp. smooth, analytic...) from the ball $B(x_0,\varepsilon)$ to a domain $D\subset\R^n$
such that
$$\psi(B(x_0,\varepsilon)\cap\Omega)=D\cap \R^n_+,\qquad
\psi(B(x_0,\varepsilon)\cap\partial\Omega)=D\cap \partial\R^n_+$$
where $\R^n_+$ is the upper half space $x_n>0$.
(If $k=0$, then diffeomorphism of class $C^0$ means homeomorphism.)

\begin{proposition}
\label{proposition-smooth}
Let $\Omega$ be a domain of class $C^0$ in the plane. Assume that Problem \eqref{eq1} has a ``classical'' solution $u$, namely: $u$ is of class $C^2$ in $\Omega$, $\Delta u=0$,
$u>0$ in $\Omega$, and
$$\forall z_0\in\partial\Omega,\quad \lim_{z\to z_0} u(z)=0\quad\mbox{ and }
\lim_{z\to z_0}|\nabla u(z)|=1.$$
Then $\Omega$ is a smooth domain (actually, real analytic). Moreover, $u$ extends
to a harmonic function defined in a neighborhood of $\Omega$.
\end{proposition}

\medskip
 
Proof: Let $z_0\in\partial\Omega$.
By the definition of a domain of class $C^0$, $z_0$ has a neighborhood
$V_{z_0}$ such that $V_{z_0}\cap\Omega$ is a Jordan domain (meaning
that its boundary is a Jordan curve).
Let $f:D^+(0,1)\to V_{z_0}\cap\Omega$ be a conformal representation on the upper half-disk. Then $f$ extends
to a homeomorphism of the closure of $D^+(0,1)$ to the closure of
$V_{z_0}\cap\Omega$ by Caratheodory's theorem
(Theorem 13.2.3 in \cite{greene}).
We may choose $f$ so that $f$ maps $(-1,1)$ to $V_{z_0}\cap\partial\Omega$.
Consider the harmonic function $\wtu(z)=u\circ f(z)$ on $D^+(0,1)$.
Then $\wtu$ is the real part of a holomorphic function $U$ on $D^+(0,1)$.
Moreover,
$$\forall x_0\in(-1,1),\quad \lim_{z\to x_0} \Re\,U(z)=0$$
so $U$ extends to a holomorphic function on $D(0,1)$ by the Schwarz reflection principle.
On the other hand, $U'=2\wtu_z$ gives
$$\frac{U'(z)}{f'(z)}=2u_z(f(z))\quad \mbox{ for $z\in D^+(0,1)$}.$$
Of course, we do not know yet that $f'$ extends continuously to $(-1,1)$.
But our hypothesis on $u$ tells us that the ratio $|\frac{U'(z)}{f'(z)}|$ does,
and morevoer,
$$\forall x_0\in(-1,1),\quad \lim_{z\to x_0}\left| \frac{U'(z)}{f'(z)}\right|=1.$$
Consequently, there exists $\varepsilon_1>0$ such that $\frac{U'(z)}{f'(z)}\neq 0$ for $z\in D^+(0,\varepsilon_1)$.
Consider the holomorphic function $h(z)=\log \frac{U'(z)}{f'(z)}$ on $D^+(0,\varepsilon_1)$. Then $\lim_{z\to x_0} \Re\,h(z)=0$ for $x_0\in(-\varepsilon_1,\varepsilon_1)$, so
$h$ extends to a holomorphic function in the disk $D(0,\varepsilon_1)$.
Hence $f'$ and $f$ extends holomorphically to the disk $D(0,\varepsilon_1)$, and $f'\neq 0$
in this disk, so $f$ is biholomorphic in a disk $D(0,\varepsilon_2)$. This implies that
the boundary of $\Omega$ is real analytic in a neighborhood of $z_0$.
Moreover, the formula $u(z)=\Re\,U(f^{-1}(z))$ shows that $u$ can be extended to
a function harmonic in a neighborhood of $z_0$. The extension is unique by
analyticity so this shows that $u$ can be extended to a neighborhood of $\Omega$.
\cqfd

\begin{remark} A similar regularity result is obtained by 
D. Khavinson, E. Lundberg and R. Teodorescu in \cite{khavinson}, Corollary
2.3, assuming that $\Omega$ is of class $C^2$.
\end{remark}

\begin{proposition}
\label{proposition-unique}
 Let $\Omega$ be an exceptional domain. Then Problem \eqref{eq1}
has a unique solution $u$.
\end{proposition}
Proof: Let $u'$ be another solution.
Then the difference $v=u-u'$ satisfies $\Delta v=0$ in $\Omega$
and $v=dv=0$ on $\partial\Omega$.
By Proposition \ref{proposition-smooth}, both $u$ and 
$u'$ extend to a neighborhood of $\Omega$.
Then the function $v_z$ is holomorphic in a neighborhood of $\Omega$ and
$v_z=0$ on $\partial\Omega$
so the zeros of $v_z$ are not isolated.
Hence $v_z\equiv 0$ and $v\equiv 0$ in $\Omega$.
\cqfd
\begin{proposition}
Let $\Omega$ be a smooth domain in the plane. Then each component of
$\partial\Omega$ is either a smooth Jordan curve or the image of a smooth
proper embedding $\gamma:\R\to\C$ (where proper means $\lim_{t\to\pm\infty}\gamma(t)=\infty$). We call the later a proper arc.
\end{proposition}
Proof. Each component of $\partial\Omega$ is a 1-dimensional
submanifold so is either diffeomorphic to the circle $\S^1$ or the real line.
In the first case, it is a smooth Jordan curve. In the later case, it is the image
of an embedding $\gamma:\R\to\C$.
We claim that $\gamma$ is proper.
If not, then
there exists a sequence $t_n\to\pm\infty$ such that $\gamma(t_n)$ has a finite limit $p$.
Now $p$ must be on the boundary of $\Omega$. The definition of a smooth domain at $p$ gives a contradiction.
\cqfd

\medskip

The condition $|\nabla u|< 1$, which we will address in Section \ref{section2},
has the following interesting geometric consequence for the domain $\Omega$:
\begin{proposition}
\label{proposition-concave}
Let $\Omega$ be an exceptional domain such that
$|\nabla u|<1$. Then
$\Omega$ is a strictly concave domain, namely: each component of $\C\setminus\Omega$ is strictly convex.
\end{proposition}
Proof. The curvature of the level set $u=0$ is given by (\cite{gray} page 72)
\begin{equation}
\label{eq-curvature-levelsets}
\kappa=\frac{u_{xx}u_y^2+u_{yy}u_x^2-2u_{xy}u_xu_y}{(u_x^2+u_y^2)^{3/2}}.
\end{equation}
Regarding sign, the curvature is positive when the curvature vector points toward $u<0$
(as can be checked in the case $u(x,y)=x^2+y^2-1$: Formula \eqref{eq-curvature-levelsets} gives $\kappa=1$).
Consider the harmonic function $g(z)=-\log |2u_z(z)|$. Then $g=0$ on $\partial\Omega$
and $g>0$ in $\Omega$.
Let $z_0\in\partial\Omega$. By rotation we may assume that $z_0=0$ and $\nabla u(0)=(1,0)$.
Then for small $\varepsilon>0$, $(0,\varepsilon)\subset\Omega$.
By the boundary maximum principle (Lemma 3.4. in \cite{GT}),
$g_x(0)>0$. On the other hand,
$$g_x(0)=-\frac{1}{|\nabla u(0)|^2}(u_x(0)u_{xx}(0)+u_y(0)u_{xy}(0))=-u_{xx}(0).$$
Hence $u_{xx}(0)<0$.
Since $u$ is harmonic we obtain $u_{yy}(0)> 0$.
Formula \eqref{eq-curvature-levelsets} gives
$\kappa(0)>0$. This means that the curvature vector of the boundary points outside of
$\Omega$, so the boundary is locally strictly concave.
Each component of $\partial\Omega$ is then globally a strictly convex curve, meaning that
it bounds a strictly convex subset of the plane.
(For the components of $\partial\Omega$ which are Jordan curves, this is standard.
For the components of $\partial\Omega$ which are proper arcs, this is also true, see
Theorem 9.40 in \cite{montiel-ros}.)
\cqfd

\medskip

Let $C$ be a convex set in the plane. For $\varepsilon>0$,
we define
$$\mu_{\varepsilon}(C)=\inf_{p\in C}\mbox{\rm Area}(C\cap D(p,\varepsilon)).$$
If $C$ is a convex set with non-empty interior then $\mu_{\varepsilon}(C)>0$.
\begin{definition}
\label{definition-thin}
Consider a subset $A$ of the plane which
is the union of a family of disjoint convex sets $(C_i)_{i\in I}$.
We say that $A$ is non-thinning if for some $\varepsilon>0$ we have
$$\inf_{i\in I} \mu_{\varepsilon}(C_i)>0.$$
\end{definition}

We will establish the correspondence in the case where the complement of
$\Omega$ is non-thinning. This prevents the components
of $\C\setminus\Omega$ from becoming  thinner and thinner.
Clearly if $\Omega$ is a concave domain with finite connectivity, or a concave periodic domain with finite connectivity in the quotient,
then its complement is non-thinning.
\section{The condition $|\nabla u|< 1$.}
\label{section2}
The goal of this section if to prove the following theorem.
\begin{theorem}
\label{theorem-finite-topology}
Let $\Omega$ be a non-trivial exceptional domain with finite connectivity (which means that $\partial\Omega$
has a finite number of components). Then
$|\nabla u|<1$ in $\Omega$.
\end{theorem}
We will in fact get more precise results: see  Theorems
\ref{theorem-bounded}, \ref{theorem-1arc} and \ref{theorem-remaining-case}.
We will also prove a result in the periodic case: see Theorem \ref{theorem-periodic}.

\subsection{The case $\partial\Omega$ compact}
We start with the case where $\partial\Omega$ is the union of a finite number
of Jordan curves.
\begin{theorem}
\label{theorem-bounded}
Let $\Omega$ be an exceptional domain such that $\C\setminus\Omega$ is bounded.
Then $|\nabla u|<1$ in $\Omega$ and
$$\lim_{|z|\to\infty} \nabla u(z)=0.$$
\end{theorem}
Proof.
Let $\wtu(z)=u(\frac{1}{z})$. Then $\wtu$ is a positive harmonic function in a punctured
disk $D^*(0,\varepsilon)$. By B\^ocher Theorem (Theorem 3.9 in \cite{axler}), we can write
$$\wtu(z)=c\log|z|+h(z)$$
where $c$ is a constant and $h$ is harmonic in the disk $D(0,\varepsilon)$.
Then
$$u(z)=-c\log|z|+h(\frac{1}{z}),$$
$$u_z(z)=-\frac{c}{2z}-h_z(\frac{1}{z})\frac{1}{z^2}.$$
Hence there exists $R$ such that $|\nabla u|<1$ for $|z|\geq R$.
Then we have $|\nabla u|\leq 1$ on the boundary of $\Omega\cap D(0,R)$.
By the maximum principle, we have $|\nabla u|\leq 1$ in $\Omega\cap D(0,R)$.
Moreover, $u_z$ is not constant in $\Omega$ so the maximum principle implies
that $|\nabla u|<1$ in $\Omega$.\cqfd

\begin{remark} D. Khavinson, E. Lundberg and R. Teodorescu (\cite{khavinson}, Theorem 3.2) have proved that under the hypothesis of Theorem \ref{theorem-bounded},
$\Omega$ is the complement of a disk, which implies the conclusion of the theorem.
\end{remark}

\subsection{A Phragmen Lindel\"of type result for univalent functions in the upper half-plane}
Next we would like to do the case where $\partial\Omega$ is a  proper arc.
For this we need a result about univalent function in the upper half-plane
$\Im\,z>0$, which we denote $H$.
(Recall that {\em univalent} means holomorphic injective.)
We have the following distortion theorem for univalent functions in the half-plane, which is an easy consequence of the standard distortion theorem for univalent functions 
in the disk:
\begin{theorem}
\label{theorem-univalent1}
Let $f$ be univalent in the upper half-plane $H$ and normalized by
$f(i)=0$ and $f'(i)=1$. Then for $z\in H$,
$$
4\frac{|z+i|-|z-i|}{(|z+i|+|z-i|)^3}\leq |f'(z)|\leq 4\frac{|z+i|+|z-i|}{(|z+i|-|z-i|)^3}.$$
\end{theorem}

Proof: Consider the M\"obius transformation
\begin{equation}
\label{eq-mobius}
w=\varphi(z)=\frac{z-i}{iz-1}
\end{equation}
This transformation is involutive and exchanges the upper half-plane $H$ with the unit disk $D(0,1)$.
Consider the function $F$ defined in the disk by $f(z)=-2F(\varphi(z))$.
Taking the derivative,
\begin{equation}
\label{eq-fprime}
f'(z)=\frac{4}{(iz-1)^2}F'(\varphi(z)).
\end{equation}
So the function $F$ is univalent in the disk and 
satisfies $F(0)=0$, $F'(0)=1$. 
The distortion theorem for normalized univalent functions in the disk (Theorem 2.5 in \cite{duren})
tells us that
\begin{equation}
\label{eq-distortion-D}
\frac{1-|w|}{(1+|w|)^3}\leq |F'(w)|\leq \frac{1+|w|}{(1-|w|)^3},\quad \mbox{ $w\in D$.}
\end{equation}
The result follows by substitution in \eqref{eq-fprime}.
\cqfd
\medskip

Using this theorem, we prove the following:
\begin{theorem}
\label{theorem-univalent2}
Let $f$ be an univalent function in the upper half-plane $H$ that extends $C^1$ to $\overline{H}$. Let $c$ be a positive number.
\begin{enumerate}
\item
If $|f'|\leq c$ on $\partial H$ then $|f'|\leq c$ in $H$.
\item If $|f'|\geq c$ on $\partial H$ then $|f'|\geq c$ in $H$.
\item If $|f'|$ is constant on $\partial H$ then $f'$ is constant in $H$.
\end{enumerate}
\end{theorem}
Proof: Without loss of generality, we may assume (replacing $f$ by $af+b$) that
$f(i)=0$ and $f'(i)=1$.
Theorem \ref{theorem-univalent1} gives
$$4\frac{|z+i|^2-|z-i|^2}{(|z+i|+|z-i|)^4}
\leq |f'(z)|\leq 4\frac{(|z+i|+|z-i|)^4}{(|z+i|^2-|z-i|^2)^3}.$$
which implies, with $z=x+iy$
\begin{equation}
\label{eq-estimate-fp}
\frac{y}{(|z|+1)^4}\leq |f'(z)|\leq \frac{(|z|+1)^4}{y^3}.
\end{equation}
The first two points of the theorem follow from Lemma \ref{lemma-univalent} below
(with $g=f'$ in the first case and $g=1/f'$ in the second one).
To prove the third point, assume that $|f'|\equiv c$ on $\partial H$.
By the first two points, $|f'|\equiv c$ in $H$. Since $f'$ is holomorphic,
it must be constant in $H$.
\cqfd
\begin{lemma}
\label{lemma-univalent}
Let $g$ be a holomorphic function on the upper half-plane $H$, continuous on $\overline{
H}$. Assume that
$$|g(z)|\leq c_1\frac{(|z|+1)^n}{y^m} \quad\mbox{ in $H$,}$$
$$|g(z)|\leq c_2 \quad \mbox{ on $\partial H$}$$
for some positive numbers $c_1,c_2$ and positive integers $n,m$.
Then $|g(z)|\leq c_2$ in $H$.
\end{lemma} 
Proof: We prove that $|g(z)|$ has polynomial growth and conclude with the
Phragmen Lindel\"of principle.
Given $R\geq 1$, consider the rectangular domain
$\Omega_R=(-R,R)\times(0,1)$ and the function
$$h_R(z)=(z-R)^m(z+R)^m.$$
We estimate the function $g(z) h_R(z)$ on $\partial\Omega_R$.
On $\{\pm R\}\times(0,1)$, we have
$$|g(z)|\leq c_1\frac{(R+2)^n}{y^m},\quad
|h_R(z)|\leq y^m(2R+1)^m\quad\Rightarrow
|g(z)h_R(z)|\leq C R^{m+n}$$
where the letter $C$ means a constant independant of $z$ and $R$.
On $[-R,R]\times\{0\}$, we have
$$|g(z)|\leq c_2,\quad |h_R(z)|\leq R^{2m}\quad \Rightarrow |g(z)h_R(z)|\leq c_2 R^{2m}.$$
On $[-R,R]\times\{1\}$, we have
$$|g(z)|\leq c_1(R+2)^n,\quad
|h_R(z)|\leq (2R+1)^{2m}\quad
\Rightarrow |g(z)h_R(z)|\leq C R^{2m+n}.
$$
Hence $|g(z)h_R(z)|\leq C R^{2m+n}$ on $\partial\Omega_R$. By the maximum principle,
$|g(z)h_R(z)|\leq C R^{2m+n}$ in $\Omega_R$.
Now if $|x|\leq \frac{R}{2}$, we have $|h_R(z)|\geq (\frac{R}{2})^{2m}$, hence
$$|g(z)|\leq C R^n.$$
This implies that in the band $\R\times[0,1]$,
\begin{equation}
\label{eq-estimate-fp2}
|g(z)|\leq C(|z|+1)^n.
\end{equation}
For $y\geq 1$, \eqref{eq-estimate-fp2} is true by the first hypothesis on $g$.
Hence \eqref{eq-estimate-fp2} holds in the whole upper half-plane.
By the Phragmen Lindel\"of principle for the half-plane (Corollary 4.2 in \cite{conway}),
$|g(z)|\leq c_2$ in $H$.
\cqfd
\subsection{The case where $\partial\Omega$ is a proper arc}
\begin{theorem}
\label{theorem-1arc}
Let $\Omega$ be an exceptional domain such that $\partial\Omega$ is a proper arc. Then $|\nabla u|=1$
in $\Omega$ and $\Omega$ is a half-plane.
\end{theorem}
Proof. Let $\gamma:\R\to\C$ be a parametrization of
$\partial\Omega$. Since $\gamma$ is proper, it extends to a continuous, injective map from the extended real line $\R\cup\{\infty\}\simeq \S^1$ to the Riemann sphere
$\C\cup\{\infty\}$. Hence $\Omega$ is bounded in the Riemann sphere by a Jordan curve.
By the Riemann mapping theorem, there exists a conformal representation $F$ from
the unit disk $D(0,1)$ to $\Omega$. By Caratheodory's Theorem
(Theorem 13.2.3 in \cite{greene}), $F$ extends
to a homeomorphism from the closed disk $\overline{D}(0,1)$ to the closure
of $\Omega$ in the Riemann sphere, namely, $\overline{\Omega}\cup\{\infty\}$.
Without loss of generality, we may assume that $F(-i)=\infty$. 
Composing with the Moebius map $\varphi$ defined in \eqref{eq-mobius},
we obtain a homeomorphism $f:\overline{H}\cup\{\infty\}\to\overline{\Omega}\cup\{\infty\}$
that is conformal in the upper half-plane $H$ and maps $\infty$ to $\infty$.
Moreover, $f(\R)=\partial\Omega$.
Since $\partial\Omega$ is smooth, $f$ extends $C^1$ to $\R$. (At this point, we know
nothing about the regularity of $f$ at infinity.)
Let $\wtu=u\circ f$. Then $\wtu$ is a positive harmonic function in $H$ with zero
boundary value. By Theorem 7.22 in \cite{axler}, $\wtu(z)=a\, \Im\,z$ for some positive constant $a$.
Then
$$\wtu_z(z)=u_z(f(z))f'(z) =\frac{-ai}{2}$$
Since $|u_z|=1$ on $\partial\Omega$, we obtain that $|f'|$ is constant on $\R$.
By Theorem \ref{theorem-univalent2}, $f'$ is constant in $H$.
 This implies that $f$ is affine, so $\Omega$ is a half-plane, and $u_z$
is constant in $\Omega$.
\cqfd
\begin{remark} F. Helein, L. Hauswirth and F. Pacard (\cite{hhp}, Proposition 6.1) have obtained this result under the additional hypothesis that $u_x>0$ in $\Omega$, a rather strong condition.
D. Khavinson, E. Lundberg and R. Teodorescu (\cite{khavinson}, Theorem 5.1) have also obtained
this result under the hypothesis that $\Omega$ is a Smirnov domain.
\end{remark}
\subsection{The remaining case}
\begin{theorem}
\label{theorem-remaining-case}
Let $\Omega$ be a non-trivial exceptional domain of finite connectivity. Assume that $\partial\Omega$ is not bounded, so contains at least one
proper arc. Then $|\nabla u|< 1$ in $\Omega$.
Moreover, for each end of $\Omega$,
$\lim \nabla u(z)$ exists and its norm is equal to one.
Finally, the number of proper arcs in $\partial\Omega$ is at most two.
\end{theorem}
Proof:
Fix some large number $R$ such that:
\begin{itemize}
\item all Jordan curves in $\partial\Omega$ are contained in $D(0,R)$,
\item all proper arcs in $\partial\Omega$ have at least one point in
$D(0,R)$,
\item $R$ is not a critical value of the function $|z|$ restricted to $\partial\Omega$.
\end{itemize}
The last point implies that $\partial\Omega$ is transverse to the circle
$|z|=R$, so intersects this circle in a finite number of points.
Consequently, $\partial\Omega\setminus D(0,R)$ has a finite number of
components. Each such component is either a curve with two endpoints on
the circle $C(0,R)$ or a proper arc with one endpoint.
Consider an unbounded component $C$ of $\Omega\setminus D(0,R)$.
By our choice of $R$, no component of $\partial C$ can be a component of
$\partial\Omega$. Hence $\partial C$ has only one component, and
we may decompose $\partial C$  as $\alpha_1\cdot\alpha_2\cdot \alpha_3$
where $\alpha_1$ and $\alpha_3$ are proper arcs $[0,\infty)\to\C$ with one endpoint
on the circle $|z|=R$, both included in $\partial\Omega$, and $\alpha_2$ is
a curve with two enpoints on $|z|=R$. (The curve $\alpha_2$ consists of arcs of the
circle $|z|=R$ together with finite parts of $\partial\Omega$).
Arguing as in the proof of Theorem \ref{theorem-1arc}, we can find a conformal
representation $f:H\to C$ which extends to a homeomorphism from
$\overline{H}\cup\{\infty\}$ to $\overline{C}\cup\{\infty\}$, mapping $\R$ to
$\partial C$ and
$\infty$ to $\infty$.

\medskip
As in the proof of Theorem \ref{theorem-1arc}, let $\wtu=u\circ f$. Then $\wtu$ is a positive harmonic function in $H$.
Since $H$ is simply connected, we may consider the conjugate harmonic
function $\wtu^*$ of $\wtu$. Consider the holomorphic function
$$U(z):=\wtu(z)+i\wtu^*(z),\quad z\in H.$$
Let $[a,b]=f^{-1}(\alpha_2)$.
Then $\Re\,U(z)=0$ on $\R\setminus[a,b]$. By the Schwarz
reflection principle, $U$ extends to a holomorphic function on $\C\setminus[a,b]$.
By the boundary maximum principle, since $\wtu>0$ in $H$, we have
$\wtu_y>0$ on $\R\setminus[a,b]$. By the Cauchy Riemann equation,
$\wtu^*_x=-\wtu_y$, so the function $\wtu^*$ is decreasing on
$(-\infty,a)$ and $(b,\infty)$. Consequently, the function $U$ takes each
pure imaginary value at most two times on $\C\setminus[a,b]$.
By Picard's theorem, $U$ has no essential singularity at $\infty$, so has at most
a pole. Since $\Re\,U$ is positive in $H$, the pole has order at most one.
This means that
we can write
$$U(z)=-ciz+h(\frac{1}{z}), \quad z\in\C\setminus [a,b]$$
where the constant $c$ is non-negative and the
function $h$ extends holomorphically at $0$.
Then
$$U'(z)=-ci -h'(\frac{1}{z})\frac{1}{z^2}.$$
If $c=0$, let $m$ be the order of the zero of $h'$ at $0$
(with $m=0$ if $h'(0)\neq 0$).
Then  there are positive constants $c_1$ and $c_2$ such that
\begin{equation}
\label{eq-U'}
\frac{c_1}{ |z|^{m+2}}\leq |U'(z)|\leq c_2\quad \mbox{ for $|z|$ large enough}.
\end{equation}
If $c>0$, then \eqref{eq-U'} still holds with $m=-2$ and $c_1=\frac{c}{2}$.
This implies in particular that $U'(z)\neq 0$ if $|z|$ is large enough.
By the Cauchy Riemann equation, we have $U'=2 \wtu_z$, so
$\wtu_z\neq 0$ for $|z|$ large enough.
From this, we conclude that for $z\in C$, $|z|$ large enough, $|\nabla u(z)|\neq 0$.
Therefore, taking a larger value of $R$ if necessary in the definition of $C$, we can assume that
$|\nabla u(z)|\neq 0$ in $\overline{C}$.
Consider the holomorphic function
$$g(z)=\frac{2\wtu_z(z)}{f'(z)}=2u_z(f(z)),\quad z\in H.$$
By \eqref{eq-estimate-fp} and \eqref{eq-U'}, we have for $|z|$ large enough
(say $|z|\geq R_0$):
$$c_1\frac{y^3}{(|z|+1)^{m+6}}\leq |g(z)|\leq c_2\frac{(|z|+1)^4}{y}.$$
Since $u_z\neq 0$ in $\overline{C}$, there exists positive constants $c'_1$ and $c'_2$ such that
for $|z|\leq R_0$,
$$c'_1\leq |g(z)|\leq c'_2.$$
Take
$$c''_1=\min\{c_1,\frac{c'_1}{R_0^3}\},\quad c''_2=\max\{c_2,c'_2R_0\}.$$
Then
$$c''_1\frac{y^3}{(|z|+1)^{m+6}}\leq |g(z)|\leq c''_2\frac{(|z|+1)^4}{y},
\qquad z\in H.$$
Since $|\nabla u|=1$ on $\alpha_1$ and $\alpha_3$, we have
$|g(z)|=1$ on $\R\setminus[a,b]$.
Since $u_z\neq 0$ on $\alpha_2$, we have by compactness that $g(z)$ and $g(z)^{-1}$ are bounded on $[a,b]$.
Hence $g(z)$ and $g(z)^{-1}$ are bounded on $\partial H$.
By Lemma \ref{lemma-univalent}, we conclude that $g(z)$ and $g(z)^{-1}$ are bounded
in $H$.
Now consider the holomorphic function
$$G(z)=\log g(z),\quad z\in H.$$
Since $g(z)$ and $g(z)^{-1}$ are bounded in $H$, $\Re\,G(z)$ is  bounded in
$H$. Moreover, $\Re\,G(z)=0$ on $\R\setminus [a,b]$.
By the Schwarz reflection principle, $G$ extends to a holomorphic function on
$\C\setminus[a,b]$.
Now $\Re\,G(z)$ is still bounded in $\C\setminus[a,b]$, so $G$ does not have
an essential singularity at $\infty$ by Picard's theorem, and cannot have a pole
either, so $G$ extends holomorphically at $\infty$.
Moreover, $G(\infty)\in i\R$.
This means that
\begin{equation}
\label{eq-limuz}\lim_{|z|\to\infty,z\in C} 2u_z(z) \mbox{ exists}
\end{equation}
and is a unitary complex number, which proves the second assertion of
Theorem \ref{theorem-remaining-case}.

\medskip

To finish the proof of the theorem,
let $C_1,\cdots,C_k$ be the unbounded components of $\Omega\setminus
D(0,R)$.
Given $\varepsilon>0$, there exists $r\geq R$ such that $|\nabla u|\leq 1+\varepsilon$
for $z\in C_i$, $|z|\geq r$.
Consider the domain
$$\Omega_r=\Omega\setminus\bigcup_{i=1}^k \{z\in C_i\,:\, |z|\geq r\}.$$
Then $\Omega_r$ is a bounded domain and $|\nabla u|\leq 1+\varepsilon$
on $\partial\Omega_r$.
By the maximum principle,
$|\nabla u|\leq 1+\varepsilon$ in $\Omega_r$.
Hence, $|\nabla u|\leq 1+\varepsilon$ in $\Omega$.
Since this holds for arbitrary positive $\varepsilon$, we have
$|\nabla u|\leq 1$ in $\Omega$. The maximum principle implies that
$|\nabla u|<1$ in $\Omega$
(else $|\nabla u|\equiv 1$ in $\Omega$ and
in this case $\Omega$ is a half-plane.)

\medskip

Since $\nabla u$ is normal to $\partial\Omega$, \eqref{eq-limuz}
implies that the normal to the
proper arc $\partial C_i$ has a limit as $t\to\pm\infty$. Hence the image of
$\partial C_i$ by $z\mapsto 1/z$ is a $C^1$ curve near zero, so the image of
$C_i$ contains a cone with vertex at the origin, positive radius $\varepsilon$ and angle
$\theta<\pi$ as close as we want to $\pi$. As the image of $\Omega$ by $z\mapsto 1/z$ can contain
at most two such cones, we have $k\leq 2$.
Theorems \ref{theorem-remaining-case} and \ref{theorem-finite-topology}
are proved.
\cqfd

\begin{remark}
D. Khavinson, E. Lundberg and R. Teodorescu (\cite{khavinson}, Corollary 4.3)
have proved the last point of Theorem \ref{theorem-remaining-case} assuming
that $\Omega$ is simply connected and Smirnov, but without requiring that $\partial\Omega$
has a finite number of components.
\end{remark}

\subsection{The periodic case}
\begin{theorem}
\label{theorem-periodic}
Let $\Omega$ be a non-trivial periodic exceptional domain invariant by
a translation $T$.
Assume that $\Omega/T$ has finite connectivity.
Then $|\nabla u|< 1$ in $\Omega$.
Moreover, the boundary of $\Omega/T$ is a finite union of Jordan curves,
$\Omega/T$ has one or two ends, each asymptotic to a half cylinder,
 and $\lim \nabla u(z)$ exists
on each end.
\end{theorem}
Proof: First observe that by uniqueness (Proposition \ref{proposition-unique}), the function $u$ is periodic:
$u\circ T=u$.
Without loss of generality we may assume that $T$ is the translation
$z\mapsto z+2\pi i$.
Let $\wtOmega=\Omega/T$.
Then $\wtOmega$ is a smooth domain in $\C/T$ (for the definition of a smooth domain is
local), so the boundary of $\wtOmega$ consists of smooth Jordan curves and proper
arcs $\gamma:\R\to\C/T$.
Choose $R>0$ large enough so that the domain $|\Re\, z|<R$ contains all the Jordan curves in
$\partial \wtOmega$.
Let $C_R$ be the half cylinder $\Re\, z>R$ in $\C/T$ and $\wtOmega_R=
\wtOmega\cap C_R$.
\begin{claim} Either
$\wtOmega_R= C_R$ or $\wtOmega_R=\emptyset$.
\end{claim}
Proof: Assume by contradiction that $\wtOmega_R$ is neither equal to $C_R$ nor
empty. Then
$(\partial\wtOmega)\cap C_R$ is not empty. By our choice of $R$, this
intersection contains a proper arc $\gamma:[0,\infty)\to C_R$ such that
$\Re\,\gamma(0)=R$. Then $C_R\setminus\gamma$ is simply connected,
so $\wtOmega_R$ lifts to an unbounded domain $\whOmega_R$ in $\C$ such that
the canonical projection $\whOmega_R\to\wtOmega_R$ is bijective.
By the proof of Theorem \ref{theorem-remaining-case}, each unbounded
component of $\whOmega_R$ contains an unbounded sector of angle
$\theta<\pi$ as close as desired to $\pi$. Since the translation
$T$ is not injective on such a sector, we get a contradiction.\cqfd

\medskip

Next assume that $\wtOmega_R$ is not empty, so is equal to $C_R$.
The function $\log:\C\setminus D(0,e^R)\to C_R$ is biholomorphic. Let $\wtu(w)=u(\log w)$, so $\wtu$ is a positive harmonic function in the domain $|w|>e^R$. By
B\^ocher theorem (Theorem 3.9 in \cite{axler}), we may write
$$\wtu(w)=c \log|w|+h(\frac{1}{w})\quad \mbox{ for $|w|>e^R$}$$
where the harmonic function $h$ extends analytically at $0$. Substituting $w=e^z$, we get
$$u(z)=c\,\Re\,z+h(e^{-z})\quad \mbox{ for $\Re\, z>R$},$$
$$u_z=\frac{c}{2}-h_z(e^{-z})e^{-z}.$$
From this, we conclude that $\lim_{x\to\infty} 2u_z(x+iy)=c$ exists.
(Note that at this point, we do not know that $|c|\leq 1$.)
This implies that $|u_z|$ is bounded in
$\wtOmega_R$.
Arguing in the same way for $x<-R$, we conclude that $|u_z|$ is bounded
in $\Omega$.
The following theorem of Fuchs \cite{fuchs} with $f=2u_z$ implies that $|2 u_z|\leq 1$ in
$\Omega$. (Indeed, the fact that $|u_z|$ is bounded rules out possibilities
(b) and (c).)
\cqfd

\begin{theorem}[Fuchs]
\label{theorem-fuchs}
Let $\Omega$ be an unbounded region of the complex plane.
If $f$ is holomorphic in $\Omega$ and $\lim\sup_{z\to\zeta,z\in\Omega} |f(z)|\leq 1$
for all $\zeta\in\partial\Omega$,
then one of the following mutually exclusive possibilities must occur:
\begin{itemize}
\item[(a)] $|f (z)|\leq 1$ for all $z\in \Omega$,
\item[(b)] $f(z)$ has a pole at $\infty$, 
\item[(c)] $\log M(r)/\log r\to\infty$ as $r\to\infty$, where
$M(r) = \sup_{|z|=r, z\in\Omega} |f(z)|$.
\end{itemize}
\end{theorem}
This is a Phragmen Lindel\"of type result. The striking fact about this result is that no
assumption is made on the geometry of the domain $\Omega$, as in the classical
Phragmen Lindel\"of principle.

\section{The minimal surface associated to an exceptional domain}
\label{section3}
To each exceptional domain $\Omega$, we associate a minimal
surface $M$ as follows.
Consider the holomorphic function $g=2 u_z$ and the holomorphic differential $dh=2 u_z dz$ on $\Omega$. Observe that both have the sames zeros, with same multiplicity.
Fix some point $z_0\in\partial\Omega$.
The Weierstrass Representation formula
\begin{equation}
\label{eq-weierstrass}
X(z)=(X_1(z),X_2(z),X_3(z))=
\Re\int_{z_0}^z \left[\frac{1}{2}(g^{-1}-g)dh,\frac{i}{2}(g^{-1}+g)dh,dh\right]
\end{equation}
defines locally a conformal, minimal immersion $X:\Omega\to\R^3$.
It turns out that $X(z)$ is in fact globally well defined in $\Omega$.
Regarding the third coordinate, we have
$$X_3(z)=\Re\int_{z_0}^z 2u_z dz=\int_{z_0}^z u_z dz+u_{\overline{z}} d\overline{z}
=\int_{z_0}^z du=u(z)$$
so $X_3(z)$ is well defined in $\Omega$.
Let
$$\psi(z)=X_1(z)+i X_2(z).$$
We will see in a moment that $\psi(z)$ is well defined in $\Omega$.
Let $M^+=X(\Omega)$. Then $M^+$ lies in the upper half-space $x_3>0$, and
the image of $\partial\Omega$ lies in the horizontal plane $x_3=0$. Since $|g|=1$
on $\partial\Omega$, we may complete $M^+$ by symmetry with respect to
the horizontal plane into a minimal surface $M$.
\begin{theorem}
\label{theorem-correspondence}
In the above setup:
\begin{enumerate}
\item $\psi(z)$ is well defined in $\Omega$.
\item For each component $\gamma$ of $\partial\Omega$, $\psi(\gamma)$ is obtained from $\gamma$ by a translation composed with conjugation. (The translation depends on the component.)
\item $M$ is a complete, immersed minimal surface in $\R^3$.
\end{enumerate}
Assume moreover that $|\nabla u|<1$ in $\Omega$. Then:
\begin{enumerate}
\setcounter{enumi}{3}
\item $\psi$ is a diffeomorphism from $\Omega$ to $\whOmega=\psi(\Omega)$ and
$M^+$ is the graph over $\whOmega$ of the function
$$\whu(z)=u(\psi^{-1}(z)).$$
Consequently, $M$ is embedded.
\end{enumerate}
Assume moreover that the complement of $\Omega$ is non-thinning
(see Definition \ref{definition-thin}). Then:
\begin{enumerate}
\setcounter{enumi}{4}
\item $\partial\whOmega=\psi(\partial\Omega)$.
\end{enumerate}
\end{theorem}
\begin{remark}
It is very much likely true that Point (5) is true without the non-thinning hypothesis
but I have not been able to prove it.
\end{remark}
Proof: A standard computation gives
\begin{equation}
\label{eq-dpsi}
d\psi=\frac{1}{2}(\overline{g^{-1}dh}-g\,dh).
\end{equation}
This gives
\begin{equation}
\label{eq-dpsi-bis}
d\psi=\frac{1}{2}(d\overline{z}-4(u_z)^2dz).
\end{equation}
We have to prove that $d\psi$ is an exact differential. In other words,
we have to prove that for any cycle $\gamma\in H_1(\Omega,\Z)$,
$\int_{\gamma} d\psi=0.$
Since $\Omega$ is a planar domain, $H_1(\Omega,\Z)$ is generated by the closed curves
in $\partial\Omega$. Let $\gamma$ be a parametrization of a component of $\partial\Omega$. Then
since $u$ is zero on $\partial\Omega$,
$$du(\gamma')=0= (u_z dz+u_{\overline{z}}d\overline{z})(\gamma').$$
Multiply by $4u_z$ and use the fact that $4u_zu_{\overline{z}}=|\nabla u|^2=1$ on $\partial\Omega$:
$$(4(u_z)^2dz+d\overline{z})(\gamma')=0$$
From this we obtain
\begin{equation}
\label{eq-dpsi2}
d\psi=d\overline{z}\mbox{ on tangent vectors to } \partial\Omega.
\end{equation}
Hence if $\gamma$ is a closed curve on $\partial\Omega$, $\int_{\gamma}d\psi=0$.
This proves Point (1). Point (2) is clearly a consequence of \eqref{eq-dpsi2}.
The metric induced on $\Omega$ by the conformal immersion $X$ is given by the standard formula
$$ds=\frac{1}{2}(|g\,dh|+|g^{-1}dh|)=\frac{1}{2}(1+|\nabla u|^2)|dz|
=\lambda(z)|dz| \quad \mbox{ with }\quad  \frac{1}{2}\leq \lambda(z)\leq 1.$$
This implies that $M$ is complete and proves Point (3).

\medskip

Proof of Point (4):
Using \eqref{eq-dpsi-bis}, the matrix of $d\psi$ is
$$\frac{1}{2}\left(\begin{array}{cc}
1-u_x^2+u_y^2 & - 2u_x u_y\\
2 u_x u_y & -1 -u_x^2+u_y^2
\end{array}\right).$$
We compute
$$\det(d\psi)=\frac{1}{4}(|\nabla u|^4 -1).$$
Since $|\nabla u|<1$, $d\psi$ is a local diffeomorphism. This implies
that the image $\whOmega=\psi(\Omega)$ is open.
The following claim proves that $\psi$ is injective, so is a diffeomorphism from
$\Omega$ to $\whOmega$.
\begin{claim}
\label{claim1}
Let $z,z'$ be two distinct points in $\Omega$. Then
$$\langle z'-z,\overline{\psi(z')}-\overline{\psi(z)}\rangle>0.$$
Here, $\langle v,v'\rangle=\Re(v\,\overline{v'})$ denotes the usual euclidean
scalar product on $\R^2$ identified with $\C$.
\end{claim}
Proof. The segment $[z,z']$ has a natural ordering which we denote $\prec$.
If $z_1,z_2$ are two points on the segment $[z,z']$ such that $z_1\prec z_2$ and
$(z_1,z_2)\subset\Omega$, then by Equation \eqref{eq-dpsi-bis}
$$\Re\left[(z'-z)(\psi(z_2)-\psi(z_1))\right]=\frac{1}{2}
\Re\left[ (z'-z)(\overline{z_2}-\overline{z_1})
-(z'-z)\int_{z_1}^{z_2}4(u_z)^2dz\right].$$
Now since $z_1,z_2$ are on the segment $[z,z']$ and $z_1\prec z_2$,
$$\Re\left[(z'-z)(\overline{z_2}-\overline{z_1})\right]
=\langle z'-z,z_2-z_1\rangle 
=|z'-z| \, |z_2-z_1|.$$
Since $|\nabla u|<1$ in $\Omega$,
$$\left|(z'-z)\int_{z_1}^{z_2} 4(u_z)^2dz\right|<|z'-z|\, |z_2-z_1|.$$
Hence
$$\langle z'-z,\overline{\psi(z_2)}-\overline{\psi(z_1)}\rangle>0.$$
If $[z,z']\subset\Omega$, the claim is proved (by taking $z_1=z$ and $z_2=z'$).
Now assume that $[z,z']$ crosses the boundary of $\Omega$.
Let $z_1$ be the first point on $[z,z']\cap \partial\Omega$
(where ``first'' refers to the ordering $\prec$ of points on $[z,z']$).
Let $\gamma_1$ be the component of $\partial\Omega$ to which $z_1$ belongs.
By Proposition \ref{proposition-concave}, $\gamma_1$ bounds a convex domain
which is in the complement of $\Omega$. The segment $[z,z']$ exits this domain
at a point $z_2\succeq z_1$ and then does not cross $\gamma_1$ anymore.
(The convexity is not crucial to this argument, but convenient).
Since $z_1$ and $z_2$ are on the same component of $\partial\Omega$,
we have by \eqref{eq-dpsi2}:
$$\langle z'-z,\overline{\psi(z_2)}-\overline{\psi(z_1)}\rangle=\langle z'-z,z_2-z_1\rangle=|z'-z|\,|z_2-z_1|.$$
Let $n$ be the number of boundary components that the segment $[z,z']$ crosses
(which must be finite by compactness).
We may find an increasing sequence of points
$z_0=z,z_1,\cdots,z_{2n},z_{2n+1}=z'$ on the segment $[z,z']$ such that
for even $i$, $(z_i,z_{i+1})\subset\Omega$ and for odd $i$,
$z_i$ and $z_{i+1}$ are on the same boundary component of $\Omega$.
By the two cases that we have seen, we have for $0\leq i\leq 2n$
$$\langle z'-z,\overline{\psi(z_{i+1})}-\overline{\psi(z_i)}\rangle >0 \quad\mbox{(unless $z_{i+1}=z_i$)}.$$
Adding all these inequalities proves Claim \ref{claim1}.

\medskip

Proof of Point (5):
Let $\whOmega=\psi(\Omega)$. Since $\psi:\Omega\to\whOmega$ is a diffeomorphism,
$\psi(\partial\Omega)\subset\partial\whOmega$ by elementary topology.
By Lemma \ref{lemma-pathologique} below (where we drop all hats),
$\partial\whOmega\subset\overline{\psi(\partial\Omega)}$.
By Claim \ref{claim-locally-compact} below, $\psi(\partial\Omega)$ is closed.
Hence $\partial\whOmega=\psi(\partial\Omega)$. 
This concludes the proof of Theorem \ref{theorem-correspondence}.\cqfd
 
 \medskip
 
\begin{claim}
\label{claim-locally-compact}
$\psi(\partial\Omega)$ is a closed subset of the plane.
\end{claim}
Note that the non-thinning hypothesis is used only to ensure that Claim \ref{claim-locally-compact} holds true.
\medskip

Proof:
let $(\gamma_i)_{i\in I}$ be the components
of $\partial\Omega$, $\widehat{\gamma}_i=\psi(\gamma_i)$,
$C_i$ the convex set bounded by $\gamma_i$ and
$\widehat{C}_i$ the convex set bounded by $\widehat{\gamma}_i$.
By Point (2) of Theorem \ref{theorem-correspondence},
$\widehat{C}_i$ is the conjugate of a translate of $C_i$.
Since $\C\setminus\Omega$ is non-thinning, there exists $\varepsilon>0$ and $\alpha>0$ such that
for all $i\in I$, $\mu_{\varepsilon}(C_i)\geq \alpha$.
Observe that $\mu_{\varepsilon}(\widehat{C}_i)=\mu_{\varepsilon}(C_i)$.
Let $z_0\in\C$. Then
$\overline{D}(z_0,\varepsilon)$ can intersect only a finite number of the convex sets
$\widehat{C}_i$, namely at most $\frac{4\pi\varepsilon^2}{\alpha}$.
(Indeed, if $p\in \widehat{C}_i\cap \overline{D}(z_0,\varepsilon)$, then $D(p,\varepsilon)\cap \widehat{C}_i$
is included in $D(z_0,2\varepsilon)$ and has area greater than $\alpha$.)
Since each $\widehat{\gamma}_i\cap \overline{D}(z_0,\varepsilon)$ is closed,
we conclude that $\psi(\partial\Omega)\cap\overline{D}(z_0,\varepsilon)$ is closed.
\cqfd
\medskip

\begin{lemma}
\label{lemma-pathologique}
Let $M$ be a complete, connected minimal surface in $\R^3$. Assume that $M$ is symmetric
with respect to the horizontal plane $x_3=0$, and that $M^+=M\cap\{x_3>0\}$ is
the graph of a function $u$ over a domain $\Omega\subset\C$. Then
$\partial\Omega\subset\overline{M^0}$, where $M^0=M\cap\{x_3=0\}$.
\end{lemma}
Of course, if $M$ is properly embedded, then $M^0$ is closed, so Lemma \ref{lemma-pathologique} says
that $\partial\Omega=M^0$. But we do not know that.
\medskip

Proof:
We follow the proof of Theorem 3.1 in \cite{daniel}.
Assume that $\partial\Omega$ contains a point $a_0$ such that $d(a_0,M^0)>0$.
Let $\varepsilon=d(a_0,M^0)$.
Choose a point $a_1\in\Omega$ such that $|a_0-a_1|\leq\frac{\varepsilon}{4}$.
 Let $a_2$ be a point in $\partial\Omega$ such that
 $|a_1-a_2|$ is minimum (which exists because $\partial\Omega$ is closed). Then
 $d(a_2,M^0)\geq\frac{\varepsilon}{2}$ and the segment $[a_1,a_2)$ is entirely
 included in $\Omega$.
 Choose a sequence of points $z_n$ on this segment
such that $z_n\to a_2$
 and $|z_n-a_2|\leq \frac{\varepsilon}{8}$.
 Let $p_n$ be the point on $M^+$ whose horizontal projection is
 $z_n$.
 Let $U_n$ be the component of $B(p_n,\frac{\varepsilon}{8})\cap M$ which
 contains $p_n$. Then for $p\in U_n$, we have
 $$d(p, M^0)\geq d(p_n,M^0)-\frac{\varepsilon}{8}
 \geq d(z_n,M^0)-\frac{\varepsilon}{8}
 \geq d(a_2,M^0)-\frac{\varepsilon}{4}
 \geq \frac{\varepsilon}{4}.$$
  Since $M^+$ is stable (as a graph), the norm of the second fundamental
 form of $U_n$ is bounded by $k=\frac{4c}{\varepsilon}$ by the estimate
 of Schoen \cite{schoen-stable}, where $c>1$ is a universal constant.
 By the uniform graph lemma (Lemma 4.1.1. in \cite{perez-ros}),
 $U_n$ is the graph over the disk $D(p_n,\frac{1}{4k})$ in the tangent plane
 $T_{p_n} M$ of a function $v_n$ which satisfies $|d^2 v_n|\leq 16 k$.
 This implies that the slope of $T_{p_n} M$ goes to infinity as $n\to\infty$, else
 the horizontal projection of $U_n$ will eventually contain $a_2$.
 Passing to a subsequence, the normal $N(p_n)$ converges to a horizontal
 vector $N_{\infty}$.
Let $\wtU_n=U_n-p_n$, so $\wtU_n$ is a minimal surface containing the point $0$.
Since it has bounded curvature, a subsequence of $\wtU_n$ converges smoothly to a minimal surface
$\wtU_{\infty}$. Moreover, the Gauss map
of $\wtU_{\infty}$ at $0$ is the horizontal vector $N_{\infty}$. I claim that $\wtU_{\infty}$
is flat. If not, then the Gauss map of $\wtU_{\infty}$ is open, so will take
values in both the upper and lower hemisphere. But then the same is true for
$\wtU_n$ for $n$ large enough, which contradicts the fact that $M$ is a graph.
Hence $\wtU_{\infty}$ is a disk of radius $\frac{1}{4k}$
in the vertical plane perpendicular to $N_{\infty}$.
This implies that the horizontal projection of $U_n$ converges to the segment
$T$ of length $\frac{1}{2k}=\frac{\varepsilon}{2c}$ centered at $a_2$ and perpendicular to $N_{\infty}$.
Then $T\subset\partial\Omega$, and since
$d(a_2,M^0)\geq \frac{\varepsilon}{2}$, we conclude that
$T\subset\partial\Omega\setminus M^0$.
The choice of $a_2$ implies that $T$ must be perpendicular to $a_2-a_1$, so
the limit normal $N_{\infty}$ is uniquely defined, up to sign.
\medskip

By changing the coordinate system, we may assume that $a_1$ and $a_2$ are on the
real axis, $a_1<0$ and $a_2=0$, so $N_{\infty}=\pm(1,0,0)$.
From what we have seen, we conclude that
for any sequence $x_n\to 0^-$, there is a subsequence such that
$\lim\frac{\nabla u(x_n,0)}{|\nabla u(x_n,0)|}=\pm (1,0)$.
Hence $|u_x(x,0)|\geq\frac{\sqrt{2}}{2}$ for $x$ close to $0$, say $x\in [-\varepsilon_1,0)$.
Consider the curve on $M$ defined by $\gamma(x)=(x,0,u(x,0))$ for $x\in[-\varepsilon_1,0)$.
Since $M$ is complete, this curve has infinite length, so
$$\infty=\int_{-\varepsilon_1}^{0}\sqrt{1+(u_x)^2}\leq 2\int_{-\varepsilon_1}^{0}|u_x|.$$
Since $u_x$ has constant sign for $x\in[-\varepsilon_1,0)$, this gives
$\lim_{x\to 0^- } u(x,0)=\pm\infty$. Since $u$ is positive, we conclude that 
the sign is $+$.

\medskip

Consider a sequence of points $z_n$ on the segment $(a_1,a_2)$ such that
$z_n\to a_2$, so $u(z_n)\to\infty$. Let $r_n=u(z_n)$.
We do the same argument again, replacing
$U_n$ by the component of $B(p_n,\frac{r_n}{2})\cap M_n$ which contains $p_n$.
Then for $p\in U_n$, we have $d(p,M^0)\geq \frac{r_n}{2}$. Fix some
arbitrary small $k>0$. By the estimate of Schoen, the norm of the fundamental
form of $U_n$ is bounded by $\frac{2c}{r_n}$ so is less than $k$ for $n$ large
enough. The above argument tells us that $\partial\Omega$ contains the segment $T$ of length $\frac{1}{2k}$ centered at $a_2$ and perpendicular to $(a_1,a_2)$.
Moreover, as each $U_n$ is a graph, $\Omega$ contains a rectangle with one side equal
to $T$ and non-empty interior (the width of this rectangle may depend on $k$).
We let $k\to 0$ and conclude that $\partial\Omega$ contains a line $L$.
By connectedness, $\Omega$ must be on one side of $L$ and $M$ is contained in
a vertical half-space.

\medskip

To prove that $u\to\infty$ on $L$,
we do the same argument again, taking $a_0$ to be any point on the line $L$.
This time we can take $a_1$ such that $(a_0,a_1)$ is perpendicular to $L$
(thanks to the existence of the above rectangle). Then $a_2=a_0$, and we
obtain that $\lim u(z)=\infty$ as $z\to L$, the limit being unifom on compact sets of $L$.
The half-space theorem of Hoffman-Meeks \cite{hoffman-meeks} gives that
$M$ is a vertical plane, which is a contradiction since $M^+$ is a graph.
(The half-space theorem of Hoffman-Meeks requires that $M$ is properly immersed.
The fact that $u\to\infty$ uniformly on compact sets of $L$ is enough, as is clear from the proof of the half-space theorem.)
\cqfd
 \section{The exceptional domain associated to a minimal bigraph}
\label{section4}
\begin{definition}
A minimal bigraph is a complete embedded minimal surface $M$ such that
$M$ is symmetric with respect to the horizontal plane $x_3=0$, and
$M^+=M\cap\{x_3>0\}$ is a graph over the domain in the horizontal plane bounded
by $M\cap\{x_3=0\}$.
\end{definition}
To each minimal bigraph $M$, we associate an exceptional domain $\Omega$ as follows.
Assume that $M^+$ is the graph of a function $\whu$ on a domain $\whOmega$.
Let $\Sigma$ be the conformal structure of $M$ (in other words, any Riemann surface
conformally equivalent to $M$).
Let $X=(X_1,X_2,X_3):\Sigma\to M$ be a conformal parametrization of $M$.
As $M$ is a minimal bigraph,
the Riemann surface
$\Sigma$ admits an antiholomorphic involution $\sigma$ corresponding
to the symmetry with respect to the horizontal plane $x_3=0$.
The fixed set of $\sigma$ divides $\Sigma$ into two components. Let
$\Sigma^+$ be the component corresponding to $M^+$.
Let $\psi=X_1+i X_2$. Then as $M^+$ is a graph over
$\whOmega$, $\psi$ is a diffeomorphism from
$\Sigma^+$ to the domain $\whOmega$.
\medskip

Let $g$ be the (stereographically projected) Gauss map of $M$ and
$dh=2\frac{\partial X_3}{\partial z}dz$ the height differential (where here $z$ denotes
a local complex coordinate on $\Sigma$).
In other words, $(\Sigma,g,dh)$ is the Weirstrass data of $M$ and $M$ is parametrized
by \eqref{eq-weierstrass}.
Assume that $M$ has been oriented so that the normal points down in $M^+$,
so that $|g|<1$ in $\Sigma^+$.
Fix some base point $p_0\in\Sigma^+$.
Define $\varphi:\Sigma^+\to\C$ by
$$\varphi(p)=\int_{p_0}^p g^{-1}dh.$$
Define $F=\varphi\circ\psi^{-1}$ and $\Omega=\varphi(\Sigma^+)$.
We have the following commutative diagram:
$$\xymatrix{
& \whOmega  \ar[dr]^{\whu} \ar[d]^{F} & \\
\Sigma^+ \ar[ru]^{\psi} \ar[r]^{\varphi} \ar@/_1pc/[rr]_{X_3} & \Omega \ar[r]^{u} & \R^{+*}\\
}
$$
% Mon premier diagramme commutatif !
\begin{theorem}
\label{theorem-inverse}
In the above setup:
\begin{enumerate}
\item $\varphi(p)$ is well defined in $\Sigma^+$.
\item For each component $\gamma$ of $\partial\whOmega$, $F(\gamma)$ is
obtained from $\gamma$ by a translation composed with conjugation.
\item $F:\whOmega\to\Omega$ is a diffeomorphism.
Moreover, for any $z,z'$
in $\whOmega$, it holds
\begin{equation}
\label{eq-enlarging}
|F(z)-F(z')|\geq |z-z'|.
\end{equation}
Consequently,
$\Omega$ is an unbounded domain whose boundary is $F(\partial\whOmega)$.
\item The function $u(z)=\whu(F^{-1}(z))$ solves Problem \eqref{eq1} in $\Omega$. Moreover, $|\nabla u|<1$ in
$\Omega$.
\end{enumerate}
\end{theorem}
Proof: 
We want to prove that $d\varphi=g^{-1}dh$ is an exact diffential on $\Sigma^+$.
Since $\Sigma^+$ is homeomorphic to a planar domain,
 it suffices to prove that $\int_{\gamma}g^{-1}dh=0$ for all closed curves
$\gamma$ on $\partial\Sigma^+$.
Let $\gamma$ be a component of $\partial\Sigma^+$. Then since $|g|=1$ on $\gamma$
and $dh(\gamma')\in i\R$, we have on $\gamma$
$$g\,dh(\gamma')=-\overline{g^{-1}dh(\gamma')}.$$
By \eqref{eq-dpsi},
\begin{equation}
\label{eq-dpsigamma}
d\psi(\gamma')=\overline{d\varphi(\gamma')}.
\end{equation}
Since $\psi$ is well defined, $d\psi$ is an exact differential, so $d\varphi$ is exact too.
This proves Point (1). Equation \eqref{eq-dpsigamma} also proves Point (2).
Regarding Point (3),
the function $g$ is holomorphic in $\Sigma^+$ and has the same zeros as $dh$ with
the same multiplicity. Hence $\varphi$ is holomorphic and $d\varphi=g^{-1}dh$ has no zero, so $\varphi$ is locally biholomorphic and $F$ is a local diffeomorphism.
We need the following
\begin{claim}
\label{claim2}
Given two distinct points $z,z'$ in $\whOmega$, we have
$$\langle \overline{F(z')}-\overline{F(z)},z'-z\rangle > |z'-z|^2.$$
Here $\langle v,v'\rangle=\Re(v\,\overline{v'})$ denotes the usual euclidean scalar product
on $\R^2$ identified with $\C$.
\end{claim}
Proof: Assume that $z_1,z_2$ are two points on the segment $[z,z']$ such that $z_1\prec z_2$
and the open segment $(z_1,z_2)$ lies inside $\whOmega$.
(Here $\prec$ denotes the natural order on the segment $[z,z']$.)
Let $\alpha:[0,1]\to\overline{\Sigma^+}$ be such that $\psi\circ\alpha$ is the constant speed parametrization of the segment $[z_1,z_2]$.
Fix some time $t\in(0,1)$ and let
$$v=\frac{1}{2}\overline{g^{-1}dh(\alpha')},\qquad
w=-\frac{1}{2}g\,dh(\alpha').$$
Then by \eqref{eq-dpsi},
$$d\psi(\alpha')=z_2-z_1=v+w,\qquad
d\varphi(\alpha')=2\overline{v}.$$
Since $|g|<1$ in $\Sigma^+$, we have $|w|< |v|$, hence
$$\langle 2v,v+w\rangle > |v+w|^2.$$
Hence
$$\langle \overline{d\varphi(\alpha')},z_2-z_1\rangle > |z_2-z_1|^2.$$
Since $z_2-z_1=\lambda (z'-z)$ with $\lambda>0$,
$$\langle\overline{d\varphi(\alpha')},z'-z\rangle > \langle z_2-z_1,z'-z\rangle.$$
Integrating from $t=0$ to $1$, we obtain
$$
 \langle \overline{F(z_2)}-\overline{F(z_1)},z'-z\rangle\ 
 =\langle\overline{\varphi(\alpha(1))}-\overline{\varphi(\alpha(0))},z'-z\rangle
 > \langle z_2-z_1,z'-z\rangle.
$$
Next assume that $z_1,z_2$ are two points on the segment $[z,z']$ such
that $z_1\prec z_2$ and $z_1$ and $z_2$ are on the same component of $\partial\whOmega$.
Then by \eqref{eq-dpsigamma}, $F(z_2)-F(z_1)=\overline{z_2}-\overline{z_1}$ so we have
$$
 \langle \overline{F(z_2)}-\overline{F(z_1)},z'-z\rangle\ = \langle z_2-z_1,z'-z\rangle.
$$
We conclude as in the proof of Claim \ref{claim1} by decomposing the segment $[z,z']$
into a finite number of segments which are either included in $\whOmega$ or
whose endpoints are on the same boundary component.
(Note that since $M$ is a minimal bigraph, the domain $\whOmega$ must
be concave.)
\cqfd

\medskip
Returning to the proof of Theorem \ref{theorem-inverse},
Claim \ref{claim2} implies Inequality \eqref{eq-enlarging}.
This implies that $F$ is injective, so $F:\whOmega\to\Omega$ is
a diffeomorphism, and that
$F$
is proper. Hence $\partial\Omega=F(\partial\whOmega)$.
Regarding Point (4), $\varphi$ is biholomorphic and $X_3$ is harmonic
so $u=X_3\circ\varphi^{-1}$ is harmonic in $\Omega$.
Since $\whu>0$ in $\whOmega$ and $\whu=0$ on $\partial\whOmega$,
we have
$u>0$ in $\Omega$ and $u=0$ on $\partial\Omega$.
Finally, differentiating $u\circ\varphi=X_3$, we get
$$2 u_z(\varphi(z))\times(g^{-1}dh)=2\frac{\partial X_3}{\partial z}dz=dh.$$
Hence
$$|\nabla u(\varphi(z))|=|g(z)|$$
which implies that $|u|<1$ in $\Omega$ and $|u|=1$ on $\partial\Omega$.
\cqfd
\section{The correspondence}
\label{section4bis}
We denote by $M[\Omega]$ the minimal surface $M$ associated to $\Omega$ by
Theorem \ref{theorem-correspondence} and by $\Omega[M]$ the domain $\Omega$
associated to the bigraph $M$ by Theorem \ref{theorem-inverse}.
Observe that the definition of $M[\Omega]$ depends on
the choice of a base point $z_0$. However, changing $z_0$ amounts to translate
$M$ by a horizontal vector. The same comment applies to $\Omega[M]$: changing
the base point $p_0$ amounts to translate $\Omega$.
Hence if we consider as equivalent two domains that differ by a translation, and
two minimal surfaces that differ by a translation, $M[\Omega]$ and
$\Omega[M]$ are well defined.
\begin{theorem}
\label{theorem-bijection}
The maps $\Omega\mapsto M[\Omega]$ and $M\mapsto\Omega[M]$ are
inverse of each other, and establish a one-to-one correspondence between
\begin{itemize}
\item exceptional domains $\Omega$ whose complement is non-thinning and such that $|\nabla u|<1$ in $\Omega$,
\item minimal surfaces $M$ which are bigraph over a domain
whose complement is non-thinning.
\end{itemize}
\end{theorem}
Proof:
\begin{itemize}
\item Assume that we are given $\Omega$ and let $M=M[\Omega]$. Recall that $M^+$ is conformally parametrized on $\Omega$ by the Weierstrass data $g=2u_z$, $dh=2u_z dz$.
Then $d\varphi=g^{-1}dh=dz$ on $\Sigma^+=\Omega$, so $\Omega[M]$
is equal to $\Omega$, up to a translation.
(Here, the conformal structure $\Sigma$ of $M$ is the ``double'' of $\Omega$,
see \cite{farkas} page 49).
\item Assume that we are given $M$ and let $\Omega=\Omega[M]$. Let $(\Sigma,g,dh)$ be the Weierstrass data
of $M$. Then $\Omega=\varphi(\Sigma^+)$ where $d\varphi=g^{-1}dh$
and
$u(\varphi(z))=X_3(z).$
By differentiating, we get
$$2u_z(\varphi(z))d\varphi=2\frac{\partial X_3}{\partial z}dz=dh,\quad\mbox{ hence } \quad
 2u_z(\varphi(z))=g(z)$$
$$\varphi^{*}(2u_z dz)=g \,d\varphi=dh.$$
Hence  $(\Sigma^+,g,dh)$ is the pullback by $\varphi$ of $(\Omega,2u_z,2u_z\,dz)$.
So $M[\Omega]=M$, up to a translation.
\end{itemize}
\cqfd

\section{Examples}
\label{section5}
In this section, we develop three examples. Please take care that
in the setup of Theorem \ref{theorem-inverse}, it is required that both $X_3>0$ and
$|g|<1$ in $\Sigma^+$. The following standard facts will be useful.
\begin{proposition}
\label{proposition-transform}
Let $(\Sigma,g,dh)$ be the Weierstrass data of a minimal surface $M$.
Then:
\begin{enumerate}
\item $(\Sigma, \frac{-1}{g},dh)$ is the Weierstrass data of $\sigma(M)$ with the
opposite orientation, where
$$\sigma(x_1,x_2,x_3)=(x_1,-x_2,x_3)$$ is the symmetry
with respect to the vertical plane $x_2=0$.
\item $(\Sigma, \frac{1+g}{1-g}, \frac{1}{2}(\frac{1}{g}-g)dh)$ is the Weierstras
data of $\rho(M)$ with the same orientation, where $$\rho(x_1,x_2,x_3)=
(-x_3,x_2,x_1)$$ is the rotation of angle $\pi/2$ around the $x_2$-axis.
\end{enumerate}
\end{proposition}
\subsection{The vertical catenoid}
\label{section-example1}
The Weierstrass data of the standard catenoid is usually written as
$$\Sigma=\C^*,\quad g=z,\quad dh=\frac{dz}{z}.$$
Then
$X_3=\log |z|$ so we see that $X_3>0$ in $|z|>1$.
Since $|g|>1$ in this domain, we use Point (1) of Proposition \ref{proposition-transform}
and take $g=\frac{-1}{z}$.
Then $\varphi(z)=-z$, so $\Omega$ is the domain $|z|>1$.

\subsection{The horizontal catenoid}
\label{section-example2}
\begin{figure}
\begin{center}
\includegraphics[height=25mm]{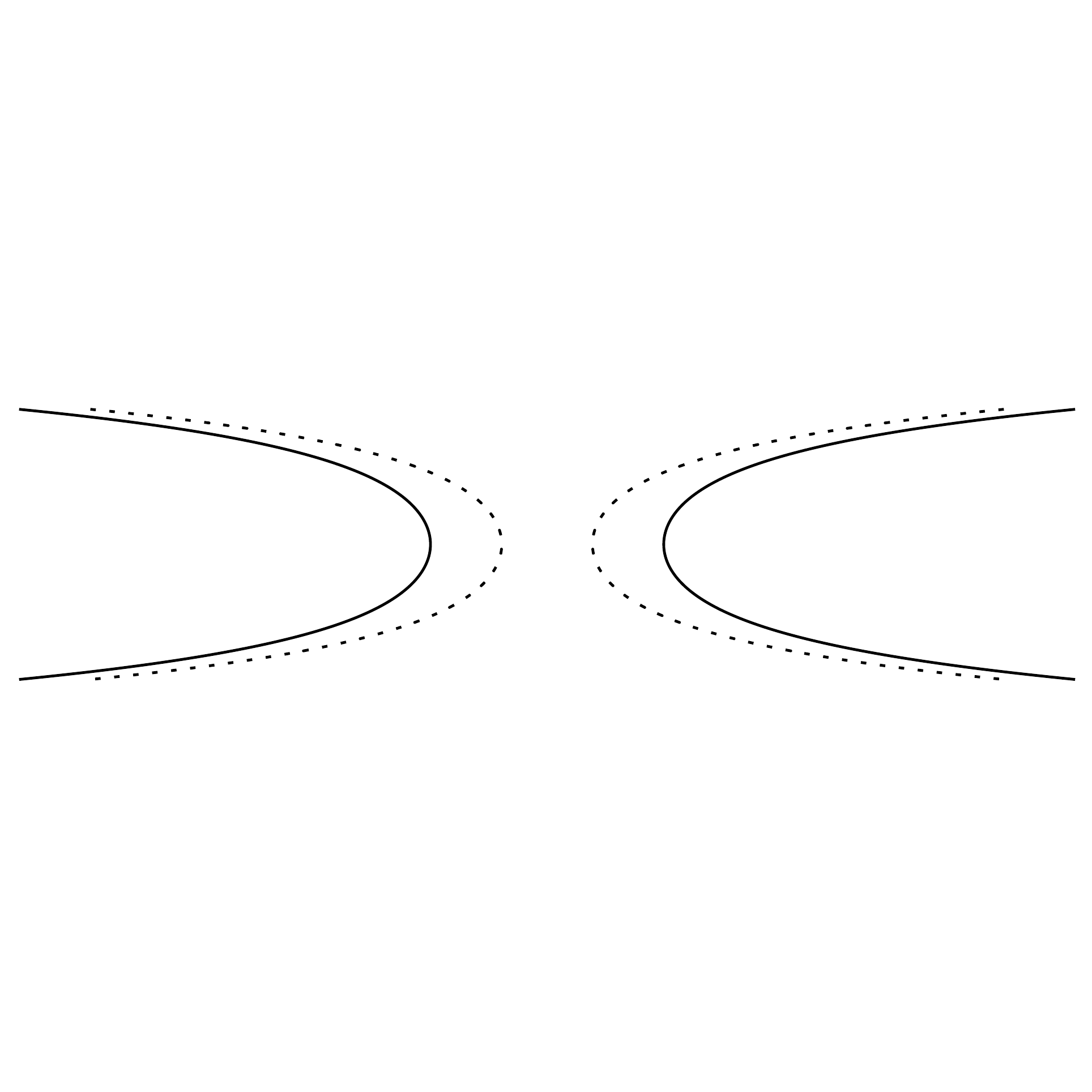}
\end{center}
\caption{The boundary of the domain $\whOmega$ over which the catenoid is a bigraph
(dots) and the boundary of the corresponding exceptional domain $\Omega$
(solid). The domains have been rotated by $90^{\circ}$ so that the figure fits the page.
Plotted with Maple.}
\label{fig1}
\end{figure}

By Point (2) of Proposition \ref{proposition-transform}, the Weierstrass data of
a horizontal catenoid is
$$\Sigma=\C^*,\qquad g=\frac{1+z}{1-z},\qquad dh=\frac{1-z^2}{2z^2}dz.$$
Here it is convenient to replace $z$ by $-z$ so
$$g=\frac{1-z}{1+z},\qquad dh=\frac{z^2-1}{2z^2}dz.$$
Then
$$X_3(z)=\frac{1}{2}\Re\left(\frac{1}{z}+z\right)$$
and 
$$X_3>0\Leftrightarrow \Re\, z>0\Leftrightarrow |g|<1.$$
$$\varphi(z)=\int g^{-1}dh=\frac{1}{2z}-\log z-\frac{z}{2}.$$
For $t$ real and $\varepsilon=\pm 1$, we have
$$\varphi(\varepsilon i e^t)=-t-\varepsilon i\left(\frac{\pi}{2}+\cosh t\right).$$
Hence, $\Omega$ is the domain $|y|<\frac{\pi}{2}+\cosh x$
(see Figure \ref{fig1}).
This is precisely the domain obtained in Proposition 2.1 of \cite{hhp}.
\begin{remark}
We see on this example that it may happen that $\varphi$ is not well defined on all of $\Sigma$.
\end{remark}
\subsection{Scherk's simply periodic surface}
\label{section-example3}
\begin{figure}
\begin{center}
\includegraphics[height=40mm]{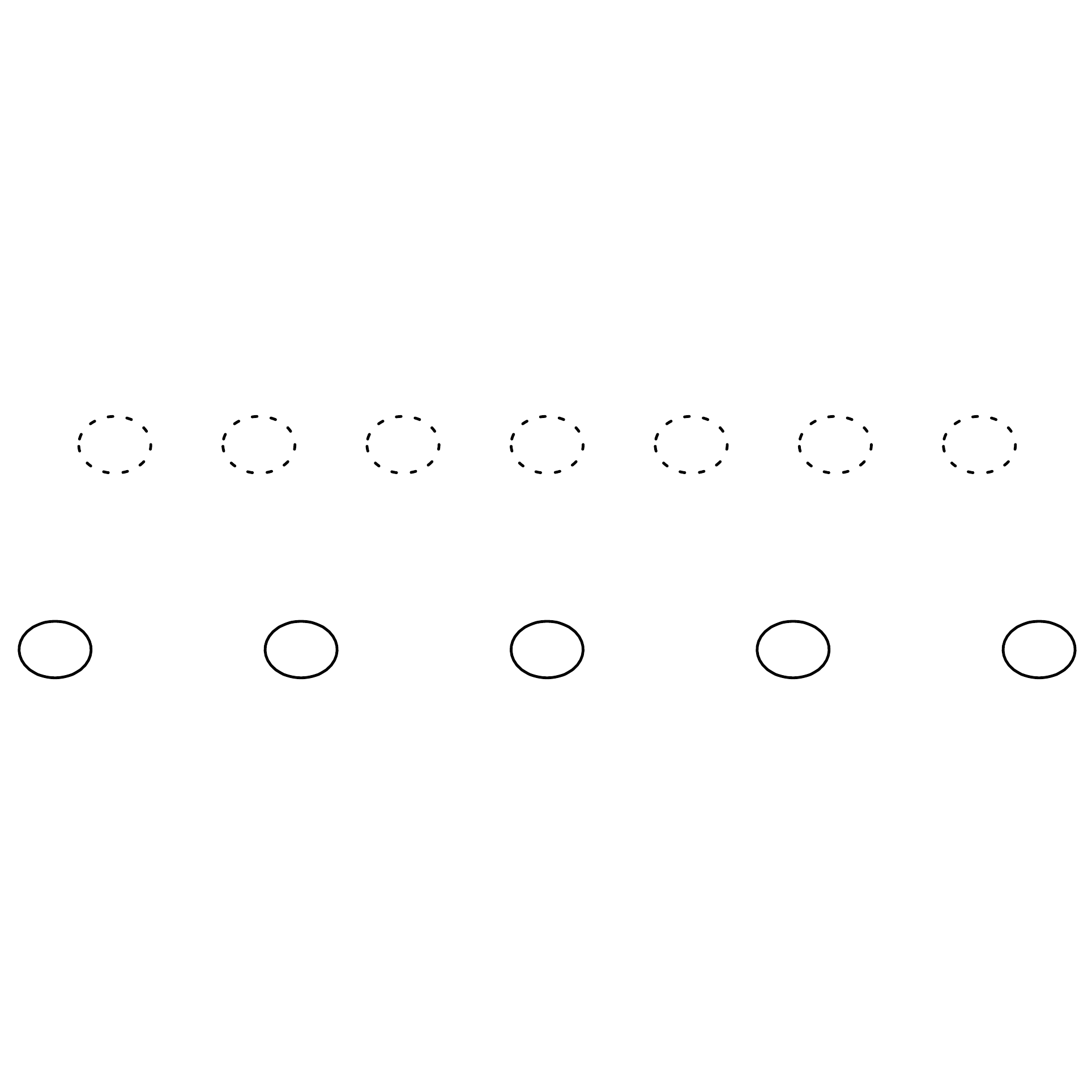}
\end{center}
\caption{Top (dots): the boundary of the domain $\whOmega$ over which the Scherk surface is a bigraph in the case $\alpha=\frac{\pi}{4}$.
Bottom (solid): the boundary of the corresponding exceptional domain $\Omega$
(translated vertically so one can see something).
Plotted with Maple.}
\label{fig2}
\end{figure}

\def\scherk{{\mbox{\tiny Scherk}}}
This periodic surface depends on a parameter $\alpha\in(0,\frac{\pi}{2})$.
Its Weierstrass data is usually written as:
$$\Sigma=\C\cup\{\infty\}\setminus\{\pm e^{i\alpha},\pm e^{-i\alpha}\},\quad
g_{\scherk}=z,\quad
dh_{\scherk}=\frac{4\sin(2\alpha)\,z\,dz}{z^4-2\cos(2\alpha) \,z^2+1}.$$
This is actually the Weierstrass data of the surface in the quotient by its period
which is the vertical vector $(0,0,2\pi)$.
The immersion \eqref{eq-weierstrass} is multi-valued on $\Sigma$ -- the multi-valuation gives rise
to the periodicity of the surface -- and is well defined on a certain covering of
$\Sigma$.
This surface is a bigraph over the vertical plane $x_1=0$ and also over $x_2=0$.
Using Point (2) of Proposition \ref{proposition-transform}, we obtain the
Weierstrass data of the horizontal Scherk surface (with horizontal
period $(2\pi,0,0)$)
which is a bigraph over the horizontal plane $x_3=0$:
$$g:=\frac{1+z}{1-z},\quad
dh=\frac{2\sin(2\alpha)\,(1-z^2)dz}{z^4-2\cos(2\alpha)\,z^2+1}.$$
Let $\sigma(z)=-\overline{z}$. Then $\sigma^*dh=-\overline{dh}$. Consequently,
taking $p_0=0$ as base point, we have $X_3(z)=0$ on $i\R\cup\{\infty\}$. From the geometry
of the Scherk surface, we know that this is precisely the zero set of $X_3$.
To determine the sign of $X_3(z)$ in $\Re\,z>0$, we observe that
$dh\simeq 2\sin(2\alpha)dz$ near $0$ so $X_3\simeq 2\sin(2\alpha)\,x$ near $0$.
Hence $X_3>0$ in $\Re\,z > 0$. Since $|g|>1$ in this domain, we use Point (1)
of Proposition \ref{proposition-transform} and replace $g$ by $-1/g$.
This gives
$$d\varphi=\frac{-2\sin(2\alpha)\,(z+1)^2dz}{z^4-2\cos(2\alpha)z^2+1}
=-\frac{(z+1)^2}{2z}dh_{\scherk}.$$
Hence
$$\Res_{e^{i\alpha}} d\varphi=-(1+\cos\alpha)\Res_{e^{i\alpha}}dh_{\scherk}=i(1+\cos\alpha).$$
The residue at $e^{-i\alpha}$ is opposite by symmetry.
Hence $\varphi$ is multi-valued on $\Sigma^+$, with multi-valuation equal to
$2\pi(1+\cos\alpha)$. So $\Omega$ is a periodic domain with period
$$T_{\alpha}=2\pi(1+\cos\alpha).$$
Now the horizontal Scherk surface is a bigraph over a domain $\whOmega$ which
is bounded by a convex curve $\gamma$ together with its translates by
multiples of $2\pi$. By Theorem \ref{theorem-inverse}, $\Omega$
is the domain bounded by $\gamma$ together with its translates by multiples of 
$T_{\alpha}$.
This is a completely explicit geometric description of $\Omega$ (see Figure \ref{fig2}).
\medskip

It turns out that one can actually compute an equation of $\gamma$.
Let us just give the main steps of the computation.
The curve $\gamma$ is the image of the circle $i\R\cup\{\infty\}$ by $\varphi$.
Then
$$\varphi(z)=-i\left[(1+\cos\alpha)\log\frac{z-e^{i\alpha}}
{z-e^{-i\alpha}}+(1-\cos\alpha)\log\frac{z+e^{i\alpha}}
{z+e^{-i\alpha}}\right].$$
Write $z=it$, $t\in\R$ and $\varphi(it)=(x(t),y(t))$. Then
$$x(t)=2\arctan\left(\frac{\sin(2\alpha)}{t^2+\cos(2\alpha)}\right),$$
$$y(t)=\cos\alpha\log\left(\frac{t^2+2\sin\alpha\,t+1}{t^2-2\sin\alpha\,t+1}\right).$$
Elimination of $t$ gives us an implicit equation of $\gamma$
\begin{equation}
\label{eq-implicit}
\cos^2\alpha\cosh(\frac{y}{\cos\alpha})=\sin^2\alpha+
\cos(2\alpha-x).\end{equation}
(More precisely, $\gamma$ is the component of the solution set of \eqref{eq-implicit}
which goes through $0$.)
In the particular case $\alpha=\frac{\pi}{4}$, Equation \eqref{eq-implicit}
simplifies to
$$\cosh(\sqrt{2}\,y)=1+2\sin x.$$
\begin{remark}
Let $\Omega_{\alpha}$ be the exceptional domain corresponding to the Scherk
surface of parameter $\alpha$. Using Equation \eqref{eq-implicit}, one can prove
that:
\begin{itemize}
\item As $\alpha\to 0$, $\frac{1}{2\alpha}\Omega_{\alpha}$ converges to the domain
$|z-1|>1$
\item As $\alpha\to\frac{\pi}{2}$, $\frac{1}{\pi-2\alpha}\Omega_{\alpha}$ converges
to the domain $|x+1+\frac{\pi}{2}|<\frac{\pi}{2}+\cosh y$.
\end{itemize}
The limit domains are, up to similitude, the examples of Sections \ref{section-example1} and \ref{section-example2}.
This corresponds to the well known fact that the horizontal Scherk surface, suitably scaled, converges to a vertical catenoid as $\alpha\to 0$ and a horizontal catenoid as $\alpha\to\frac{\pi}{2}$.
As a consequence, if we consider as equivalent two domains which differ by a similitude, we
can put all the examples of Sections \ref{section-example1}, \ref{section-example2} and 
\ref{section-example3} in a continuous family $\Omega_{\alpha}$
for $\alpha\in[0,\frac{\pi}{2}]$: $\Omega_0$ is the example of Section \ref{section-example1}
and $\Omega_{\frac{\pi}{2}}$ is the example of Section \ref{section-example2}.
\end{remark}

\section{Classification results}
\label{section6}
\begin{theorem}
\label{theorem-classification1}
Let $\Omega$ be an exceptional planar domain of finite connectivity (meaning that
$\partial\Omega$ has a finite number of components).
Then $\Omega$ is one of the following domains:
\begin{itemize}
\item a half-plane,
\item the outside of a disk,
\item the domain $|y|<\frac{\pi}{2}+\cosh x$, 
up to a similitude.
\end{itemize}
\end{theorem}
Assume that $\Omega$ is not a half-plane.
By Theorem \ref{theorem-finite-topology}, $|\nabla u|< 1$ in $\Omega$.
Let $M$ be the minimal bigraph associated to $\Omega$ by Theorem
\ref{theorem-correspondence}.
If $\partial\Omega$ is compact, then $M$ has two ends.
By Theorem \ref{theorem-bounded}, $\lim_{z\to\infty} g(z)=0$. This implies that
$M$ has finite total curvature.  
By a theorem of Schoen \cite{schoen}, an embedded minimal surface with finite total
curvature and two ends is a catenoid. Since $g=0$ at the top end, $M$ is a
vertical catenoid.
\medskip

Else, let $k\geq 1$ be the number of proper
arcs in $\partial\Omega$. Then $M$ has $k$ ends.
By Theorem \ref{theorem-remaining-case}, $k\leq 2$ and
the limit of the Gauss map at each end exists
and is a complex number of norm 1. This implies that $M$ has finite total
curvature.
If $k=1$, then $M$ is vertical plane because the only embedded minimal surface with
finite total curvature and one end is the plane. This is not possible because
the vertical plane is not a bigraph.
Hence $k=2$, and $M$ is a horizontal catenoid by the theorem of Schoen.
\cqfd

\begin{theorem}
\label{theorem-classification2}
Let $\Omega$ be a periodic exceptional domain.
Assume that $\Omega$ has finite connectivity in the quotient.
Then $\Omega$ is one of the following domains:
\begin{itemize}
\item a half-plane,
\item the exceptional domain corresponding to a horizontal Scherk surface
(namely, one of the domains described in Section \ref{section-example3}, up to similitude).
\end{itemize}
\end{theorem}
Proof:
Assume that $\Omega$ is not a half-plane.
By Theorem \ref{theorem-periodic}, $|\nabla u|<1$ in $\Omega$. Let
$M$ be the minimal bigraph associated to $\Omega$ by Theorem
\ref{theorem-correspondence}. Then $M$ is a
periodic minimal surface with horizontal period $T$. By
Theorem \ref{theorem-periodic}, $M^+/T$ is bounded by a finite number of
Jordan curves in the plane $x_3=0$, and has at most two ends.
Moreover, the Gauss map has a limit at each end, so $M/T$ has finite total curvature.
By a theorem of Meeks-Rosenberg \cite{meeks-rosenberg}, 
the ends of $M/T$ are either of planar,
helicoidal or Scherk type.
In both the planar and helicoidal cases, $M^+/T$ would intersect the horizontal plane $x_3=0$ in a non-compact set. Hence $M/T$ has at most four Scherk-type ends.
If $M/T$ has two Scherk-type ends then it is a plane, which is not possible. So it has four Scherk-type ends.
By a theorem of Meeks-Wolf \cite{meeks-wolf}, $M$ is a Scherk surface.\cqfd
\begin{remark}
The theorem of Meeks-Wolf is a difficult result. Moreover, Theorem \ref{theorem-classification2} is equivalent to the Theorem of Meeks-Wolf: Indeed, using the Alexandrov
moving plane method, one can prove that a periodic minimal surface with 4 Scherk-type ends must be a minimal bigraph over some plane. For this reason, I don't think that
there is an elementary proof of Theorem \ref{theorem-classification2}.
\end{remark}
\section{Immersed domains}
\label{section7}
In \cite{hhp}, the authors also propose to study Problem \eqref{eq1} on arbitrary flat Riemannian manifolds with boundary.
They construct examples which have some analogy with immersed minimal surfaces
called $k$-noids. This was another hint at the correspondence between exceptional
domains and minimal surfaces.
The correspondence, however, does not generalize to arbitrary flat surfaces. Let
me propose a setup where the correspondence extends. This will allow us to recover
the examples discussed in \cite{hhp}, and more.
The following definitions are standard:
\begin{definition}
\label{definition-multidomain}
\begin{enumerate}
\item A (smooth, 2-dimensional) immersed domain $\Omega$ is a smooth, complete, flat, 2-dimensional Riemannian manifold-with-boundary such that there exists a map
$f:\Omega\to\C$ which is a local  isometry, called the
developing map of $\Omega$.
\item We say that $\Omega$ has embedded ends if the developing map is injective on each end of $\Omega$.
\end{enumerate}
\end{definition}
Note that by definition of a
manifold-with-boundary, $\Omega$ includes its boundary. We will denote by
$\mathring{\Omega}=\Omega\setminus\partial\Omega$ the set of interior points of $\Omega$.
Here the word {\em complete} means that $\Omega$ is complete as a metric space.
A flat Riemannian manifold always admits locally a developing map,
but the developing map is in general not globally defined unless the manifold
is simply connected. The definition of an immersed domain requires the developing map
to be globally defined.
\medskip

An immersed domain $\Omega$ (with non-empty boundary) is called
exceptional if Problem \eqref{eq1} has a solution
$u$ on $\Omega$ (where $\Delta u$ and $|\nabla u|$ are
computed for the metric of $\Omega$).
Theorem \ref{theorem-finite-topology} generalizes to:
\begin{theorem}
Let $\Omega$ be a non-trivial exceptional immersed domain with finite connectivity and
embedded ends. Then $|\nabla u|<1$ in $\mathring{\Omega}$.
\end{theorem}
Proof: Theorem \ref{theorem-finite-topology} is proved by showing that
$|\nabla u|$ is bounded in each unbounded component of $\Omega\setminus
D(0,R)$. Since we assume that our immersed domain has embedded ends, the
proof carries over.\cqfd

\medskip

Next we recall the definition of {\em strong symmetry} from \cite{cosin-ros}, Definition 1.
Let $X:M\to\R^3$ be an isometric immersion of a connected orientable surface $M$,
and $\Pi$ be a plane in $\R^3$ which we normalize as the horizontal plane $x_3=0$.
Denote by $S$ the symmetry with respect to the plane $x_3=0$, and
$$M^+=M\cap\{X_3>0\},\quad
M^-=M\cap\{X_3<0\},\quad
M^0=M\cap\{X_3=0\}.$$
\begin{definition}
$M$ is strongly symmetric with respect to $\Pi$ if:
\begin{enumerate}
\item There exists an isometric involution $s:M\to M$ such that
$X\circ s=S\circ X$,
\item The set of fixed points of $s$ is $M^0$,
\item The third coordinate $N_3$ of the Gauss map of $M$ takes positive (resp.
negative) values on $M^+$ (resp. $M^-$).
\end{enumerate}
\end{definition}
With these definitions, Theorem \ref{theorem-bijection} generalizes to:
\begin{theorem}
\label{theorem-rnoids}
There is a one-to-one correspondence between the following
two classes of objets:
\begin{itemize}
\item immersed domains $\Omega$ which have finite connectivity,
embedded ends and are homeomorphic to a planar domain,
\item complete, immersed minimal surfaces $M$ which are strongly
symmetric, have finite total curvature, embedded ends, and such that
$M^+$ is homeomorphic to a planar domain.
\end{itemize}
\end{theorem}
There are plenty of such minimal surfaces.
The basic example is the Jorge-Meeks $k$-noid, which has $k\geq 3$ horizontal
catenoidal ends.
Genus zero examples with $k\geq 3$ horizontal catenoidal ends are classified
by C. Cosin and A. Ros in \cite{cosin-ros}, they form a $2k-2$ parameters family
which includes the Jorge-Meeks $k$-noid as the most symmetric member.
The corresponding exceptional domains are the domains constructed in Section
4 of \cite{hhp}.
Genus one examples with $k\geq 3$ horizontal catenoidal ends are constructed
by L. Mazet in \cite{mazet}.
Pictures of higher genus examples can be seen on the minimal surface archive
of M. Weber \cite{weber}.
\medskip

Proof:
Assume that we are given an exceptional immersed domain $\Omega$ with developing
map $f$, satisfying all the hypothesis of Theorem \ref{theorem-rnoids}.
We define the holomorphic differential $dh$ by
$$dh=2 d^{(1,0)} u=2u_z dz$$
where $z$ is a local conformal coordinate on $\Omega$.
The holomorphic function $g$ is defined by
$$g=\frac{dh}{df}=\frac{2u_z}{f'}.$$
Since $f$ is a local isometry, the metric of $\Omega$ in the local coordinate $z$
is $|f'(z)dz|$.
Hence
$$|\nabla u|=\frac{|2u_z|}{|f'|}=|g|.$$
Let $M^+$ be the minimal surface parametrized on $\Omega$ by the Weierstrass
Representation formula \eqref{eq-weierstrass}.
To see that the immersion $X$ is well defined, consider the differential
$$d\psi=dX_1+i \,dX_2=\frac{1}{2}(\overline{g^{-1}dh}-gdh).$$
Then on the boundary of $\Omega$ we have
$$d\psi=\overline{g^{-1}dh}=\overline{df}.$$
Since the developing map is well defined in $\Omega$, $df$ is an exact differential.
Since $\Omega$ is homeomorphic to a planar domain, $H_1(\Omega,\Z)$ is
generated by the closed curves in $\partial\Omega$. Hence $d\psi$ is an
exact differential on $\Omega$ and $X$ is well defined.
Since $X_3=0$ and $|g|=1$ on $\partial\Omega$, we may extend $M^+$
by symmetry with respect to the plane $x_3=0$ into a strongly symmetric immersed
minimal surface $M$. The metric induced on $\Omega$ by the conformal immersion $X$
is
\begin{equation}
\label{eq-ds}
ds=\frac{1}{2}(|g^{-1}dh|+|gdh|)=|df|\frac{1+|g|^2}{2}.
\end{equation}
Since $|g|<1$ in $\Omega$, this implies that $M$ is complete.
Let $E$ be an end of $\Omega$. There are two cases:
\begin{itemize}
\item If $f(E)$ is the complement of a bounded domain in $\C$, then $E$ is conformally
a punctured disk and $df$ has a double pole at the puncture.
Moreover, by Theorem \ref{theorem-bounded}, $|\nabla u|\to 0$ so $g$ has a zero
at the puncture. This implies that $g^{-1}dh$, $gdh$ and $dh$ have at most double
poles at the puncture. Since this characterizes embedded ends of finite total curvature,
we conclude that the corresponding end of $M$ is embedded.
\item If the boundary of $f(E)$ is not bounded: then by passing to a sub-end, we
may assume that $f(E)$ is a concave domain bounded by
$\alpha_1\cdot\alpha_2\cdot\alpha_3$, where $\alpha_1$ and $\alpha_3$ are proper arcs: $[0,\infty)\to\C$
and $\alpha_2$ is a straight segment connecting the endpoints of $\alpha_1$ and $\alpha_3$. The proof of Point (4) of Theorem \ref{theorem-correspondence} says that
$\psi$ is injective on $E$, so $X(E)$ is a graph and the corresponding end of $M$ is
embedded.
 (Indeed, if $z,z'$ are two points in $f(E)$, the segment $[z,z']$
can only cross the boundary components $\alpha_1$ and $\alpha_3$.)
\end{itemize}
Finally, the fact that $g$ has a limit at each end implies that $M$ has finite total curvature.

\medskip

Conversely, assume that we are given a minimal surface $M$ satisfying all the hypothesis
of Theorem \ref{theorem-rnoids}.
Let $(\Sigma,g,dh)$ be the Weierstrass data of $M$. Since $M$ is strongly symmetric,
the Riemann surface $\Sigma$ admits a antiholomorphic involution $s$ such that
$X\circ s=S\circ X$. Moreover, the fixed set of $s$ divides $\Sigma$ into two
components $\Sigma^+$ and $\Sigma^-$ such that $|g|<1$ in $\Sigma^+$ and
$|g|>1$ in $\Sigma^-$.
(Observe that since $M$ has finite total curvature, it is properly immersed, so the strong halfspace theorem of Hoffman Meeks (Theorem 2 in \cite{hoffman-meeks}) implies that
$M^+$ is connected.)
Consider the differential $df=g^{-1}dh$ in $\Sigma^+$.
Since $g$ and $dh$ have the same zeros with same multiplicity, $df$ is holomorphic
with no zeros in $\Sigma^+$.
On the boundary of $\Sigma^+$, we have $df=\overline{d\psi}=dX_1-i\, dX_2$. Since $X$
is well defined and $\Sigma^+$ is homeomorphic to a planar domain,
this implies that $df$ is exact. By integration, we obtain a well defined
holomorphic function $f:\Sigma^+\to \C$ with non-zero derivative.
We define $\Omega$ as $\Sigma^+$ with the conformal metric $|df|$ and $f$ as developing map.
Formula \eqref{eq-ds} shows that the metric $|df|$ is equivalent to the metric $ds$
induced by the immersion $X$ on $\Sigma^+$. Hence
$\Omega$ is complete, so is an immersed domain.

\medskip

It remains to prove that $\Omega$ has embedded ends. Fix an end $E$ of $M$.
Since $M$ is complete and has finite total curvature, Osserman's theorem tells us that
$E$ can be parametrized conformally on a punctured disk. Moreover, $g$ and $dh$ extend meromorphically at the puncture.
Since $M$ is strongly symmetric,
Lemma 4 in \cite{cosin-ros} tells us that the asymptotic normal at the end is either
horizontal or vertical.
Therefore, either $g=0$, $g=\infty$ or $|g|=1$ at the puncture corresponding to the end.
\begin{itemize}
\item If $g=0$ at the puncture, then the end can be parametrized on a punctured disk entirely
included in $\Sigma^+$. Moreover, since the end is embedded, $g^{-1}dh$ has a double
pole at the puncture. Therefore, $f$ has a simple pole, so is injective in a neighborhood
of the end. This implies that the corresponding end of $\Omega$ is embedded (and asymptotic to a plane).
\item If $g=\infty$ at the puncture, then the end can be parametrized on a punctured disk
entirely contained in $\Sigma^-$, so we do not see it in $\Omega$.
\item If $|g|=1$ at the puncture, then Point (b) of Lemma 4 in \cite{cosin-ros}
says that the end is asymptotic to a horizontal catenoid. Therefore,
$E^+=E\cap M^+$ is a graph over a concave domain $\widehat{E}$ in the plane, which we may take to be bounded by $\alpha_1\cdot\alpha_2\cdot\alpha_3$, where $\alpha_1$ and
$\alpha_3$ are convex curves included in $M^0$ and $\alpha_2$ is a straight segment.
The proof of Point (3) of Theorem \ref{theorem-inverse} tells us that $f$ is injective
in $E^+$.
\cqfd
\end{itemize}

\bigskip
\noindent
{\sc Martin Traizet\\
Laboratoire de Math\'ematiques et Physique Th\'eorique\\
Universit\'e Fran\c cois Rabelais\\
37200 Tours, France.\\}
{\em email address: }\verb$martin.traizet@lmpt.univ-tours.fr$

\end{document}